%
%

\magnification=1200

\font\titfont=cmr10 at 12 pt

\font \fr = eufm10


\font\AAA=cmr14 at 12pt
\font\BBB=cmr14 at 8pt

\overfullrule=0in

  \def\harr#1#2{\ \smash{\mathop{\hbox to .3in{\rightarrowfill}}\limits^{\scriptstyle#1}_{\scriptstyle#2}}\ }

\def\half{\hbox{${1\over 2}$}}
\def\smfrac#1#2{\hbox{${#1\over #2}$}}
\def\oa#1{\overrightarrow #1}
\def\dim{{\rm dim}}
\def\dist{{\rm dist}}

\def\deg{{\rm deg}}

\def\log{{\rm log}}
\def\Hess{{\rm Hess}}
\def\Sym{\rm Sym}

\def\tr{{\rm tr}}
\def\max{{\rm max}}

\def\span{{\rm span\,}}
\def\Hom{{\rm Hom\,}}

\def\End{{\rm End}}

\def\arr{\longrightarrow}

\def\fpsh{{\rm PSH}(X,\f)}
\def\Core{{\rm Core}}
\def\dis{f_M}


\def\Theorem#1{\medskip\noindent {\AAA T\BBB HEOREM \rm #1.}}
\def\Prop#1{\medskip\noindent {\AAA P\BBB ROPOSITION \rm  #1.}}
\def\Cor#1{\medskip\noindent {\AAA C\BBB OROLLARY \rm #1.}}
\def\Lemma#1{\medskip\noindent {\AAA L\BBB EMMA \rm  #1.}}
\def\Remark#1{\medskip\noindent {\AAA R\BBB EMARK \rm  #1.}}
\def\Note#1{\medskip\noindent {\AAA N\BBB OTE \rm  #1.}}
\def\Def#1{\medskip\noindent {\AAA D\BBB EFINITION \rm  #1.}}

\def\Ex#1{\medskip\noindent {\AAA E\BBB XAMPLE \rm    #1.}}
\def\Qu#1{\medskip\noindent {\AAA Q\BBB UESTION \rm    #1.}}

\def\pf{\medskip\noindent {\bf Proof.}\ }
\def\qed{\hfill  $\vrule width5pt height5pt depth0pt$}

\def\df{d^{\phi}}
\def\hk{\_{\rm l}\,}
\def\n{\nabla}
\def\w{\wedge}

   \def\cc{{\cal C}}     
   \def\cp{{\cal P}}
   
\def\ce{{\cal E}}   
\def\ch{{\cal H}}   
   \def\cn{{\cal N}}
\def\cd{{\cal D}}
\def\cl{{\cal L}}
\def\cp{{\cal P}}
\def\cf{{\cal F}}

\def\gerG{{\fr{\hbox{g}}}}

\def\vf{\varphi}

\def\wt{\widetilde}
\def\wh{\widehat}

\def\and{\qquad {\rm and} \qquad}
\def\arr{\longrightarrow}
\def\ol{\overline}
\def\bbr{{\bf R}}\def\bbh{{\bf H}}\def\bbo{{\bf O}}
\def\bbc{{\bf C}}

\def\bbz{{\bf Z}}

\def\a{\alpha}
\def\b{\beta}
\def\d{\delta}
\def\e{\epsilon}
\def\f{\phi}

\def\l{\lambda}
\def\o{\omega}

\def\s{\sigma}
\def\x{\xi}

\def\D{\Delta}
\def\L{\Lambda}
\def\G{\Gamma}
\def\O{\Omega}

\def\fp{$\phi$-plurisubharmonic }
\def\fh{$\phi$-pluriharmonic }
\def\pfp{ partially \fh}

\def\psh{plurisubharmonic}
\def\lloc{L^1_{\rm loc}}

\def\bo{\partial \Omega}
\def\fc{$\phi$-convex }
\def\PSH{\rm PSH}

\def\BM{\lambda}
\def\Der{D}
\def\CH{{\cal H}}

\def\fdim{{\rm fd}(\f)}
\def\rn{{\bbr}^n}

\def\B{C}
\def\YY{0}
\def\AA{1}
\def\ZZ{2}
\def\KK{3}
\def\BB{4}
\def\CC{5}
 
\def\FF{6}


\ 
\vskip .3in

\centerline{\titfont AN INTRODUCTION TO POTENTIAL THEORY   }
\smallskip

\centerline{\titfont IN CALIBRATED GEOMETRY }
\bigskip

\centerline{\titfont F. Reese Harvey and H. Blaine Lawson, Jr.$^*$}
\vglue .9cm
\smallbreak\footnote{}{ $ {} \sp{ *}{\rm Partially}$  supported by
the N.S.F. }

\vskip .3in
\centerline{\bf ABSTRACT} \medskip
  \font\abstractfont=cmr10 at 10 pt

{{\parindent= .7in\narrower\abstractfont \noindent
In this paper we introduce and study the notion of plurisubharmonic functions
in  calibrated geometry.  These functions generalize 
the classical plurisubharmonic functions from complex geometry
 and enjoy  their important properties. 
 Moreover, they exist in abundance whereas the corresponding
 pluriharmonics are generally quite scarce.
 A number of the results
established in complex analysis via plurisubharmonic functions are extended
to calibrated manifolds.  This paper introduces and investigates  questions
of  pseudo-convexity in the context of a general calibrated manifold $(X,\f)$.
Analogues of totally real submanifolds are introduced and used to construct
enormous families of strictly $\f$-convex
spaces  with every topological type allowed by 
Morse Theory. Specific calibrations are used as examples throughout.

 ${\!\!\!\!\!\!\!\!\!\!\!\!\!\!\!\!\! {\rm In \ a \ sequel}}$,   the duality between \fp functions and $\f$-positive
 currents is investigated.  
 This study involves  boundaries, 
generalized Jensen measures, and other geometric objects on a calibrated manifold.

}}

\vfill\eject\
\vskip 1in

\centerline{\bf TABLE OF CONTENTS} \bigskip

\qquad\qquad\qquad\qquad  \YY. Introduction.\smallskip

\qquad\qquad\qquad\qquad  \AA. Plurisubharmonic Functions.\smallskip

\qquad\qquad\qquad\qquad  \ZZ.  The $\f$-Hessian. \smallskip

\qquad\qquad\qquad\qquad  \KK. Elliptic Calibrations.\smallskip

\qquad\qquad\qquad\qquad  \BB.  Convexity in Calibrated Geometries.\smallskip

\qquad \qquad\qquad\qquad  \CC. Boundary Convexity.\smallskip

\qquad\qquad\qquad\qquad  \FF. $\f$-Free Submanifolds.\bigskip

\centerline{\bf Appendices}
\medskip

\qquad\qquad\qquad\qquad  Appendix A: Submanifolds which are $\f$-Critical.  \smallskip

\qquad\qquad\qquad\qquad  Appendix B: Constructing $\f$-Plurisubharmonic Functions.
\smallskip

\qquad\qquad \qquad\qquad Appendix C: Structure of the Core.
 \vfill\eject

\centerline{\bf \YY. Introduction.}
\medskip

Calibrated geometries, as introduced in [HL$_1$],  are  geometries
of distinguished submanifolds  determined by a fixed, closed differential form $\f$ on
a riemannian manifold $X$.
The basic example is that of a K\"ahler manifold (or more generally a symplectic manifold,
with compatible almost  complex structure) where the distinguished submanifolds are the  holomorphic curves.
However, there exist many other interesting  geometries, each carrying a wealth of $\f$-submanifolds, particularly  on spaces with special holonomy.  The relationship between spinors and calibrations established in [DH] provides additional interest.  Calibrated manifolds  have attracted particular attention in recent years due to their appearance in generalized Donaldson theories ([DT], [Ti]) and in modern versions of string theory in Physics ([GLW], [GP], [AFS], [Her], [MS], [G], [EM], [GW], [MC] for example).

Unfortunately, analysis on these spaces $(X,\f)$ has been difficult,  in part because  
there is generally no reasonable analogue of
the  holomorphic functions  and  transformations which exist in the  K\"ahler  case.
However, in complex analysis there are many important results which can be established
using only the plurisubharmonic  functions (cf. [Ho], [D]).   It turns out that analogues of these functions 
exist in abundance on any calibrated manifold, and they enjoy almost all the pleasant properties 
of their cousins from complex analysis. The point of this paper is to introduce and study 
these functions and related notions of convexity. 

In a sequel   [HL$_2$]   these notions
will be related to $\f$-positive currents and their boundaries, 
generalized Jensen measures, and other geometric objects on a calibrated manifold.

We begin by defining our notion of    $\f$-plurisubharmonicity for smooth functions
on any calibrated manifold $(X,\f)$.  In the K\"ahler case  they
are exactly the classical plurisubharmonic functions.  We then study the basic properties of these functions, and subsequently use them to establish a series of results in geometry 
and analysis on $(X,\f)$.  

A  fundamental result is that:\smallskip

\centerline{\sl The restriction of a \fp function to a $\f$-submanifold $M$}
\centerline{\sl is subharmonic in the induced metric on $M$.}\smallskip

Any  convex function on the riemannnian manifold $X$ is  $\f$-plurisubharmonic. Moreover, at least locally, there exists an abundance 
of \fp  functions which are not convex. 

The definition of  $\f$-plurisubharmonicity extends from smooth functions
to arbitrary distributions on $X$.
Such distributions enjoy
all the nice properties of generalized subharmonic functions.
In this paper, however,  we shall focus mainly on the smooth case, except for Section 3.

To define \fp functions on a calibrated manifold $(X,\f)$ where deg$(\f)=p$, we introduce a second order
differential operator
$
\ch^\f : C^\infty(X) \ \to\ \ce^p(X),
$
the {\sl $\f$-Hessian}, given by  
$$
\ch^\f(f)\ =\ \BM_\f(\Hess f)
$$
where $\Hess f$ is the riemannian hessian of $f$ and $\BM_\f:\End(TX) \to \L^pT^*X$
is the bundle map given by  $\BM_\f(A) = D_{A^*}(\f)$ where  $D_{A^*}:\L^pT^*X\arr \L^pT^*X$ is the natural extension of $A^*:T^*X \to T^*X$ as a derivation.

When the calibration $\f$ is parallel there is a natural factorization
$$
\ch^\f \ =\ d d^\f
$$
where $d$ is the de Rham differential and $d^\f: C^\infty(X) \ \to\ \ce^{p-1}(X)$ is given by
$$
d^\f f\ \equiv \ \nabla f \hk \f.
$$
In general these operators are related by the equation: $\ch^\f  f = dd^\f f - \nabla_{\nabla f} (\f)$.

Recall that a calibration $\f$ of degree $p$  is a closed $p$-form  with the property that $\f(\x)\leq 1$ for all unit simple tangent $p$-vectors $\x$ on $X$. Those $\x$ for which $\f(\x)=1$ are called 
$\f$-planes, and the set of $\f$-planes is denoted by $G(\f)$. With this understood, a function
$f\in C^\infty(X)$ is defined to be {\bf  \fp} if $\ch^\f(f)(\x) \geq 0$ for all $\x\in G(\f)$. It is {\bf strictly
\fp } at a point $x\in X$ if $\ch^\f(f)(\x) > 0$ for all $\f$-planes $\x$ at $x$.  In a similar fashion, $f$ 
is called {\bf $\f$-pluriharmonic} if $\ch^\f(f)(\x) = 0$ for all $\x\in G(\f)$.
Denote by $\PSH(X,\f)$ the convex cone of \fp functions on $X$.

When $X$ is a complex manifold with a K\"ahler form $\o$, one easily computes that
$d^\o = d^c$, the conjugate differential.  In this case, $\ch^\o = dd^\o =dd^c$ and the 
$\o$-planes correspond to the complex lines in $TX$. Hence, the definitions above coincide
with the classical notions of plurisubharmonic and pluriharmonic functions on $X$.

With this said, we must remark that in many calibrated manifolds the $\f$-pluriharmonic
functions are scarce. For the  calibrations on manifolds with strict G$_2$ or Spin$_7$ holonomy,  for example,  every pluriharmonic function is constant.  For the Special Lagrangian
calibration $\f = {\rm Re}\{dz\}$, every $\f$-pluriharmonic  function $f$ defined locally in $\bbc^n$ is of the  form $f=a+q$ where $a$ is affine and $q$ is  a traceless Hermitian quadratic form
(cf.  [Fu] and Proposition 1.10.) Nevertheless, as we stated above, the \fp functions in any calibrated geometry are locally abundant.

The fundamental property of the $\f$-Hessian: 
$$
\left(\ch^\f f\right)(\x)\ =\ {\rm trace}\left\{ \Hess f\bigr|_{\x}  \right\} \qquad \ \ {\rm for \ all \ } \f-{\rm planes\ } \x
$$
is established in Section \ZZ \ (Corollary \ZZ.5).  This gives the useful fact that
\medskip
\centerline {   $f$ is $\f$-plurisubharmonic \ \ 
$\iff \ \ \tr_\xi \left\{ \Hess f\bigr|_{\x}  \right\} \ \geq\ 0 \quad \forall \, \x\in G(\f)$.}
\medskip

\medskip
\centerline
{\AAA E\BBB LLIPTIC  \AAA C\BBB ALIBRATIONS.}\smallskip

This brief section is an introduction to the theory of \fp distributions.  A very mild condition
on the calibration is needed to ensure ``ellipticity'', namely, $G(\f)$ should involve all the variables.
(This is stronger than requiring that the calibration $\f$ involve all the variables. See Example 2 in Section 3.)
Under this assumption, each \fp  distribution is, in fact, $\lloc$ (locally Lebesgue integrable) and has a canonical point-wise representative which is $[-\infty, \infty)$-valued and upper semi-continuous,
given by the limit of the means over balls.  The usual properties of plurisubharmonic function in complex analysis are valid for these \fp functions.  See [HL$_4$] for a more comprehensive development,
which is also calibration independent.

\medskip

Beginning with Section \BB\  the \fp functions are used to study geometry and analysis 
on calibrated manifolds.  The first concept to be addressed is the analogue
of pseudoconvexity in complex geometry.

\vfill\eject
\centerline
{\AAA C\BBB ONVEXITY.}\smallskip

 Let  $(X,\f)$ be a calibrated manifold and $K\subset X$ a closed
subset. By the {\bf $\f$-convex hull of $K$ } we mean the subset
$$
\wh K\ =\ \{x\in X : f(x) \leq \sup_K f \ \ {\rm for \ all \  } f\in \PSH(X,\f)\}
$$
The manifold $(X,\f)$ is said to be {\bf  $\f$-convex} if $K\subset \subset X \ \Rightarrow
\  \wh K \subset \subset X$  for all $K$.

\Theorem {\BB.3} {\sl A calibrated manifold  $(X,\f)$ is $\f$-convex if and only if 
it admits a \fp proper exhaustion function $f:X\to \bbr$.} \medskip

The manifold $(X,\f)$ will be called {\bf strictly $\f$-convex} if it admits an exhaustion function $f$ which is  strictly $\f$-plurisubharmonic, 
and it will be called  {\bf strictly $\f$-convex at infinity}  if $f$ is strictly \fp  outside of a compact subset.
It is shown that in the second case,  $f$ can be assumed to be
\fp everywhere. 
Analogues of Theorem \BB.3 are established in each of these cases.

Note that in complex geometry, strictly \fc manifolds are Stein and manifolds which are strictly \fc at infinity are called strongly pseudoconvex.

We next consider the {\bf core} of $X$ which is defined to be the set of points
$x\in X$ with the property that no $f\in\PSH(X,\f)$ is strictly \fp at $x$. The following results are established:
\smallskip

1)\ The manifold $X$ is strictly $\f$-convex at infinity if and only if Core$(X)$ is compact.
\smallskip

2)\ The manifold $X$ is strictly $\f$-convex if and only if Core$(X)=\emptyset$.
\smallskip

Examples of complete calibrated manifolds with compact  cores 
are given in the final subsection of \S \BB. 
A very general construction of strictly \fc manifolds is presented in \S \FF.
We next examine the analogues of pseudoconvex boundaries in calibrated geometry.
\bigskip
\centerline
{\AAA B\BBB OUNDARY \AAA C\BBB ONVEXITY.}\smallskip

Let $\O\subset X$ be a domain with smooth boundary $\bo$, and let $\rho:X\to \bbr$ be a {\bf defining function for $\partial \O$}, that is, a smooth function defined on a neighborhood of $\overline \O$ with $\O=\{x : \rho(x) <0\}$, and $\nabla \rho \neq 0$ on $\bo$. Then $\bo$ is said
to be {\bf $\f$-convex} if 
$$
\ch^\f(\rho)(\x)\ \geq \ 0 \qquad {\rm for\  all\ } \f-{\rm planes\  } \x {\rm \ tangential \ to\ }  \bo,
\eqno{(\YY.1)}
$$
i.e., for all $\x\in G(\f)$ with span$(\x)\subset T(\bo)$. The boundary $\f$ is {\bf strictly $\f$-convex} if the inequality in (\YY.1) is  strict everywhere on $\bo$.  These conditions are independent of the
choice of defining function $\rho$.

\Theorem{\CC.6} {\sl Let $\O\subset\subset X$ be a compact domain with strictly $\f$-convex boundary, and let $\delta= -\rho$ where $\rho$ is an arbitrary defining function for $\bo$.
Then $-\log \delta:\O\to \bbr$ is strictly \fp outside a compact subset.  In particular,
the domain $\O$ is strictly $\f$-convex at infinity.
}
\medskip

Elementary examples show that the converse of this theorem does not hold in general.
However, there is  a weak partial converse.
\Prop{\CC.9}  {\sl Let $\O\subset\subset X$ be a compact domain with smooth boundary.
Suppose $\f$ is parallel and consider the function $\d = {\rm dist}(\bullet, \bo)$.
If $-\log \d$ is \fp near $\bo$, then $\bo$ is $\f$-convex.}

\medskip

We note that boundary convexity can be  interpreted geometrically as follows.
Let $II$ denote the second fundamental form of the hypersurface $\bo$ oriented by the 
outward-pointing normal.  Then $\bo$ is \fc if and only if trace$(II\bigr|_{\x})\leq 0$
for all $\f$-planes $\x$ which are tangent to $\bo$.
 In the strict case one also has the following.

\Theorem {\CC.14}  {\sl  Let $(X,\f)$ be a strictly \fc manifold and $\O\subset\subset X$ a 
domain with smooth boundary.  Then the following are equivalent.
\smallskip

(i)  \ \ $\bo$ is strictly \fc.
\smallskip

(ii) \ $\tr_\x\left\{ II_{\bo}\right\} <0$ for all tangential $\f$-planes  $\x$.
\smallskip

(iii)   There exists a smooth defining function $\rho$ for $\bo$ which is strictly
\fp

\qquad  on a neighborhood of $\overline\O$.   }

 \bigskip
 
\centerline
{\AAA $\f$-F\BBB REE \AAA S\BBB UBMANIFOLDS  \BBB AND \AAA S\BBB TRICTLY
$\f$-\AAA C\BBB ONVEX \AAA S\BBB UBDOMAINS.}\smallskip
\medskip

We next examine the analogues in calibrated geometry of the totally real
submanifolds in complex analysis.  Using the methods of [HW$_{1,2}$] 
we then show how to construct strictly \fc manifolds in enormous families with every
topological type allowed by Morse theory.

Let $(X,\f)$ be any fixed calibrated manifold.
A closed submanifold $M\subset X$ is called {\bf $\f$-free} if there are no $\f$-planes
tangential to $M$, i.e.,    no $\x\in G(\f)$ with $\span \x \subset TM$.

Note that $M$ is automatically $\f$-free if it is {\bf $\f$-isotropic}, that is, if $\f\bigr|_M\equiv 0$. . 

Any submanifold of dimension $< p$ is $\f$-free, and generic local submanifolds
of dimension $p$ are $\f$-free. Furthermore, any submanifold of a $\f$-free submanifold is again $\f$-free.   

The {\bf $\f$-free dimension} of $(X,\f)$, denoted $\fdim$, is defined to be the largest dimension of a 
$\f$-free vector  subspace of $T_xX$ for $x\in X$.  The first result is the following
generalization of the Andreotti-Frenkel Theorem [AF] for Stein manifolds.

\Theorem {\FF.2}  {\sl   Suppose $(X,\f)$ is a strictly $\f$-convex manifold. Then $X$ has the homotopy type of a CW complex of dimension $\leq \fdim$.}
\medskip

For a K\"ahler manifold $(X,\o)$ of complex dimension n,  the $\o$-free dimension  is $n$ and the     
$\o$-free submanifolds are those which are { totally real}  (e.g., the Lagrangian submanifolds).
Furthermore the $\o^p/ p!$-free dimension is $n+p-1$ and a 
submanifold $M$ is  $\o^p/ p!$-free if there are no complex $p$-planes tangent to $M$ at any 
point.

In Special Lagrangian geometry on an $n$-dimensional  Calabi-Yau manifold $(X,\o,\f)$,  the $\f$-free dimension is $2n-2$ 
and the $\f$-free submanifolds  are exactly the symplectic submanifolds (e.g., the complex hypersurfaces).

 For a quaternionic K\"ahler manifold $(X,\Psi)$ 
 of dimension 4$n$, where $\Psi = {1\over 6} (\o_I^2 + \o_J^2 + \o_K^2 )$
 is the fundamental 4-form,  the $\Psi$-free dimension is 3$n$.  
 For the higher degree calibrations $\Psi_p \equiv {1\over (2p+1)!}  (\o_I^2 + \o_J^2 + \o_K^2 )^p$
 the free dimension is $3(n-p+1)$.
 
If $(X,\Phi)$ is an 8-dimensional Spin$_7$-manifold with Cayley calibration $\Phi$,
${\rm fd}(\Phi)=4$.

If $(X,\f)$ is a 7-dimensional G$_2$-manifold with associative calibration $\phi$, then
 ${\rm fd}(\phi)=4$.   So if $X$ is \fc, it has homotopy dimension $\leq 4$.
Recently, I.Unal [U] has shown that {\sl for every connected manifold $M$ of dimension $< 4$
(compact or non-compact) there exists  a strictly \fc G$_2$-manifold which is 
homotopy-equivalent to $M$.}

The relationship between $\f$-free submanifolds and convexity is expressed
in the next two results.

\Theorem {\FF.4}  {\sl Suppose $M$ is a closed submanifold of $(X,\f)$ and let 
${\rm dist}_M^2(x)$ denote   the square of the distance to $M$.  Then
$M$ is $\f$-free  if and only if 
the function
${\rm dist}_M^2(x)$ is strictly \fp at each point in $M$ (and hence in a neighborhood of $M$).
}\medskip

The existence of $\f$-free   submanifolds ensures the existence of many strictly \fc domains in $(X,\f)$.

\Theorem{\FF.6}  {\sl  Suppose $M$ is a $\f$-free submanifold of $(X,\f)$.
  Then there exists a fundamental   system  $\cf(M)$ of strictly $\f$-convex neighborhoods of $M$, 
  each of which admits a  deformation retraction onto $M$.
  }

\medskip

This result provides rich families of strictly convex domains.  
The neighborhoods in $\cf(M)$ include the  sets  
$\{x: {\rm dist}_M(x)<\e(x)\}$ for positive functions $\e\in C^\infty(M)$ which die arbitrarily rapidly at infinity.
As noted,   any submanifold of  dimension $<p$ is $\f$-free. Furthermore, 
any submanifold of a $\f$-free submanifold is again $\f$-free. 

For example if $X$ is a Calabi-Yau manifold with Special Lagrangian calibration $\f$,
then any symplectic submanifold $Y\subset X$ is $\f$-free, as is any smooth submanifold
$A\subset Y$. The topological type of such manifolds $A$ can be quite complicated.

This construction can be refined even further by replacing the submanifold $A\subset Y$ with an arbitrary closed
subset.
It turns out that  the following two classes of subsets:
\smallskip

\qquad (1)\ Closed subsets $A$ of $\f$-free submanifolds
\smallskip

\qquad (2)\ Zero sets  of non-negative strictly \fp\ functions $f$
\smallskip

\noindent
are essentially  the same.

 \medskip

We mention that the operator $d^\f$ has been independently found by M. Verbitsky 
[V] who studied the generalized K\"ahler theory (in the sense of Chern) on G$_2$-manifolds.
The authors would like to thank Robert Bryant for useful comments and conversations
related to this paper.


\vfill\eject


\centerline{\bf \AA. Plurisubharmonic Functions}
\medskip

Suppose $\phi$ is a calibration on a riemannian manifold $X$.  
The  $\phi$-Grassmannian, denoted $G(\phi)$, consists of the unit simple
$p$-vectors $\xi$ with $\phi(\xi)=1$, i.e., the $\phi$-planes.  An oriented submanifold $M$ is  a
{\sl $\phi$-submanifold}, or is {\sl calibrated by $\phi$}, if the
oriented unit tangent space ${\oa T}_x M$ lies in $G_x(\phi)$ for each
$x\in M$, or equivalently, if $\phi$ restricts to $M$ to be the volume form
on $M$.  Let $n = \dim X$ and $p = {\rm degree}(\phi)$. 

\Def {\AA.1} The {\bf $\df$-operator} is defined by $$\df f \ =\ \nabla f
\hk      \phi$$ for all smooth functions $f$ on $X$. 

\medskip
\noindent
Hence
$$
\df : \ce^0(X)\ \arr\ \ce^{p-1}(X) \and d \df :\ce^0(X)\ \arr\ \ce^{p}(X)
$$
where   $\ce^p(X)$ denotes the space of $C^\infty$ $p$-forms on $X$. 
This $d d^\f$ operator  provides a way of defining
plurisubharmonic functions in   calibrated geometry when the calibration $\f$ is parallel.

If  $\o$ is a K\"ahler form on a complex manifold, then
$
d^\o \ =\ d^c\ =\ -J\circ d
$
is the conjugate differential. Thus, the $dd^\f$-operator generalizes the 
$dd^c$-operator in complex geometry.
Although no analogue of  a holomorphic function exists on a calibrated
manifold, there is an analogue of the real part  of a holomorphic function.

\Def {\AA.2} Suppose $\nabla \phi=0$. A  function $f\in C^\infty(X)$
is {\bf $\phi$-plurisubharmonic} if 
$$
(d \df f)(\xi)\ \geq 0 \qquad {\rm for\  all } \ \ \xi \in G(\phi).
$$
The set of such functions will be denoted $\PSH(X,\phi)$.
If $(d \df f)(\xi)>0$ for all $\xi \in G(\phi)$, then $f$ is {\bf strictly  $\phi$-plurisubharmonic}. 
If $(d \df f)(\xi)=0$ for all $\xi \in G(\phi)$, then $f$ is  {\bf $\f$-pluriharmonic}.
Finally, $f$ is {\bf partially $\f$-pluriharmonic} if $f$ is \fp and, at each point, there
exists a $\f$-plane $\x$ with $(dd^\f f)(\x)=0$.    
\medskip
Note that $f$ is \fh if and only if both $f$ and $-f$ are $\f$-plurisubharmonic,  and that $f$ is \pfp if and only if 
$f$ is \fp but not strict at any point.

\Remark {\AA.3}  If $\f$ is not parallel, we define \fp functions by replacing $dd^\f f$,
in Definition \AA.2, with  $$\ch^\f(f) = dd^\f f - \nabla_{\nabla f} \phi.$$  This modified $dd^\f$-operator is discussed in detail in Section \ZZ. Note that the difference $\nabla_{\nabla f} \f$ is a first order operator.

\medskip

\Ex{}   Consider the Special Lagrangian calibration
$\f={\rm Re}( dz)$ on $\bbc^n$.  Let $Z_{ij}$ denote the bidegree $(n-1,1)$ form obtained
from $dz = dz_1\wedge\cdots\wedge dz_n$ by replacing $dz_i$ with $d{\bar z}_j$  (in the $i$th position).    
A short calculation shows that
$$
dd^\f f \ =\ 2 {\rm Re} \biggl\{   \sum_{i,j=1}^n{\partial^2 f \over \partial {\bar z}_i \partial {\bar z}_j}
Z_{ij}\biggr\}
+(\D f) {\rm Re}( dz)
\eqno{(\AA.1)}
$$

\medskip\noindent
 {\AAA N\BBB OTE THAT:} 
\  1) The constant functions are \fh.
\smallskip

\qquad\qquad 2) \ If $a,b>0$ and $f,g\in\fpsh$, then $af+bg\in \fpsh$.
\medskip

The next result justifies the use of the word plurisubharmonic in the
context of  a $\phi$-geometry.  A calibration $\f$  is {\bf integrable} if for
each point $x\in X$ and each $\x\in G_x(\f)$ there exists a $\f$-submanifold $M$ through $x$ with 
$\oa{T_xM} =\x$.

\Theorem{\AA.4} {\sl  Let $(X,\f)$ be any calibrated manifold.  If a function  $f\in C^\infty(X)$ is $\phi$-plurisubharmonic, then the restriction of $f$ to any $\phi$-submanifold $M\subset X$ is subharmonic in the induced  metric.  If $\f$ is integrable, then the converse holds.}
\medskip

Theorem \AA.4  is an immediate consequence of the formula
$$
\ch^\f f\bigr|_M\ =\  \left(    dd^\f f - \nabla_{\nabla f} \f   \right)\bigr|_M\ =\ (\D_M f )\, {\rm vol}_M
\eqno{(\AA.2)}
$$
This formula follows  from the three equations (\ZZ.7), (\ZZ.12), and (\ZZ.15), proved
below, and the fact that $\f$-submanifolds are minimal submanifolds.  We continue for the moment to present results whose proofs will be given  in  Section \ZZ.

The $\phi$-plurisubharmonic functions enjoy many of the
useful properties of their classical cousins in complex analysis.
 The next result is useful, in particular if one wishes to only consider smooth \fp functions

\Lemma {\AA.5} {\sl Let $f,g\in C^\infty(X)$ be
$\phi$-plurisubharmonic.\smallskip

a)\ \  If $\psi\in C^\infty(\bbr)$ is convex and increasing, then $\psi \circ
f$ is $\phi$-plurisubharmonic.\smallskip

b) \ \ If $\psi\in C^\infty(\bbr^2)$ is convex, and is  increasing in each variable, then $\psi(f,g)$

\qquad is \fp.}
\pf See Appendix B.
\medskip

\Remark{\AA.6}  
Part b) can be used to construct a \fp smoothing $h_\e$ of the maximum $h=\max\{f,g\}$ of two 
\fp functions $f,g$ with:
\smallskip

1)\ \ $h_\e$ decreasing as $\e\to 0$, \smallskip

2)\ \ $h_\e  -\e\leq \max\{f,g\} \leq h_\e$  for all $\e>0$, \smallskip

3)\ \ $h_\e  = \max\{f,g\} $ on the set where $|f-g|\geq \e$.
\medskip

\noindent
To see this, note first that $\max\{t_1,t_2\} = \half(t_1+t_2)+\half|t_1-t_2|$. 
Now choose   a convex function $\vf \in C^\infty(\bbr)$  with $\vf(0)=\half$, $|\vf'|\leq 1$,  and
$\vf(s)\geq |s|$  with equality when  
$ |s|\geq 1$.  Then $\vf_\e(s)=\e\vf({s\over \e})$ provides a smooth approximation to the function
$|s|$, namely $\vf_\e(s) -\e \leq |s| \leq \vf_\e(s)$.
  The function $\psi_\e(t_1,t_2) = \half(t_1+t_2)+\half\vf_\e(t_1-t_2)$ approximates
$\max\{t_1,t_2\}$, and the function $h_\e= \psi_\e(f,g)$ approximates $h=\max\{f,g\}$.  To complete the proof, note that \smallskip
$
\qquad\qquad
{\partial \psi_\e\over \partial t_1}\ =\ {1\over2}\left(1+\vf'\left({t_1-t_2\over \e}\right)\right),
\qquad{\rm  } \qquad
{\partial \psi_\e\over \partial t_2}\ =\ {1\over2}\left(1-\vf'\left(  {t_1-t_2\over \e}\right)\right)
$
\smallskip
\noindent
and

\centerline{
$
2\e\Hess \psi_\e\ =\ \vf''\left({t_1-t_2\over \e}\right)\left(\matrix{1&-1\cr-1&1\cr}\right).
$
}

\medskip

While pluriharmonic functions are often scarce,  the partially pluriharmonic  functions represent the calibrated analogue of solutions to the homogeneous
Monge-Amp\'ere equation, and they are sufficiently abundant to solve the Dirichlet
Problem [HL$_{3,6}$].

For now we mention a ``fundamental''  example.

\Prop{\AA.7} {\sl Suppose $\f \in \L^p\bbr^n$ is a parallel calibration.  Set 
$$
E(x) \ = \ \log |x|\ \  {\rm if\ } p=2\and
E(x) \ = \ -{1\over( p-2)}  \, {1\over |x|^{p-2}} \ \  {\rm if\ } p\geq3
\eqno{(\AA.3)}
$$
\noindent
Then $E$ is \fp on $\bbr^n-\{0\}$.  Moreover,  $E$ is \pfp on $\bbr^n-\{0\}$ if and only if
 each unit vector $e\in \bbr^n$ is contained in a $\f$-plane $\x\in G(\f)$.}

\Remark{{\bf (The Abundance of \fp Functions)}}   We shall see in the next section that
any convex function on the riemannian manifold $X$ is automatically $\f$-plurisubharmonic.
However, there always exist  huge families of locally defined \fp functions which are not convex.
This follows, for example,   from duality considerations as in Remark \ZZ.9 below.
However, in section $\FF$ we give a general construction of \fp functions from any $\f$-free submanifold, which shows that such functions exist in abundance.



\bigskip

\centerline{\bf    Pluriharmonic Functions} 

\medskip
The \fh functions are a natural replacement for the holomorphic functions in complex geometry.
 However, while  \fp functions 
are abundant, the $\f$-pluriharmonic functions are often quite scarce.  To illustrate this 
phenomenon we shall sketch some of the basic facts in the ``classical'' cases.

To begin we note that for some calibrations $\f$, one has that:
$$
dd^\f f\ =\ 0 \qquad {\rm if\ and\ only\ if\ }\qquad (dd^\f f)(\x) \ =\ 0 \quad {\rm for\ all\  } \x\in G(\f)
\eqno{(\AA.4)}
$$
while for others this is not true. It is the right hand side that defines pluriharmonicity.
If (\AA.4) holds and the basic map $\BM_\f$, defined in section \ZZ, is everywhere injective (as in Example \AA.14), then the only pluriharmonic functions are the affine functions, i.e., the functions with parallel gradient.
Note that if $f$ is affine, then $\nabla f$ splits the manifold locally as a riemannian 
product $X=\bbr\times X_0$.  

\Ex{\AA.8. {\bf (Complex geometry)}}  Let $\o$ be a K\"ahler form on a complex manifold
$X$.  Then
$
d^\o = d^c
$
is the conjugate differential, $dd^cf$ is the complex hermitian Hessian of $f$, $G(\o)$ is the grassmannian of complex lines, and the statement (\AA.4) is valid. In particular, the 
$\o$-pluriharmonic functions are just the classical pluriharmonic functions on $X$.

For the higher divided powers $\Omega_p = {1\over p!}\o^p$ one computes that 
$d d^{\Omega_p} f = \Omega_{p-1} dd^c f$. Furthermore, it can be deduced  from the discussion
in Remark \ZZ.13 that (\AA.4) holds in this case. Therefore, the $\Omega_p$-pluriharmonic functions
are also just the classical pluriharmonic functions.

\Ex{\AA.9. {\bf (Quaternionic-K\"ahler  geometry)}}  
Let $\bbh$ denote the quaternions and consider $\bbh^n$ as a right-$\bbh$
vector space.  Each of the complex structures $I, J, K$ (right multiplication by $i,j,k$) determines a
 K\"ahler form $\o_I, \o_J, \o_K$ respectively.  The 4-form
 $$
 \Psi\ \equiv\  {1\over6}(\o_I^2 + \o_J^2 + \o_K^2)
 \eqno{(\AA.5)}
 $$
on $\bbh^n\equiv \bbr^{4n}$ is a calibration with $G( \Psi)$ consisting of the oriented quaternion
lines in $\bbh^n$.  In this case, $d d^ \Psi f \equiv 0$ if and only if $\Hess f \equiv 0$.  However,
the assertion (\AA.4) is not valid in this case, and in fact there is a rich family of $ \Psi$-pluriharmonic functions.  For example, if $f$ is $\o_I$-pluriharmonic,  then $f$ is $ \Psi$-pluriharmonic.
Hence, so is any $\o$-pluriharmonic $f$ where $\o=a\o_I+b\o_J+c\o_K$ with $a^2+b^2+c^2=1$.

It is well known that the only $ \Psi$-submanifolds in $\bbh^n$ are the
affine quaternion lines.

Of course the calibration (\AA.5) exists on any quaternionic K\"ahler manifold,
i.e.,  one with ${\rm Sp}_n \cdot {\rm Sp}_1$-holonomy.  (See [GL] for examples.)   With this 
full  holonomy group it seems unlikely that there are many
$\Psi$-pluriharmonic functions. However, if the holonomy is contained in ${\rm Sp}_n$, they
exist in abundance as seen in the next example.

\Ex{\AA.10. {\bf (Hyper-K\"ahler  manifolds)}}  Let $(X, \o_I,\o_J,\o_K)$ be a hyper-K\"ahler manifold.
Then $X$ carries  several parallel calibrations. There are, of course, the K\"ahler forms
$\o=a\o_I+b\o_J+c\o_K$ with $a^2+b^2+c^2=1$, and two others of particular interest.

(1) \ Let $\Psi = {1\over6}(\o_I^2 + \o_J^2 + \o_K^2)$.  Then as in Example \AA.9, {\sl any 
$\o$-pluriharmonic function is $\Psi$-pluriharmonic}.  Hence, the sheaf of $\Psi$-pluriharmonic
functions is quite rich on any manifold with SP$_n$-holonomy.
 On the other hand there are precious few $\Psi$-submanifolds.

(2)\  Consider the generalized Cayley form 
$
\Xi \ \equiv\ {1\over2}(\o_I^2-\o_J^2-\o_K^2).
$
For this calibration there exist no interesting pluriharmonic functions, at least in dimension 8, but there
are many $\Xi$-submanifolds (cf. [BH]).

\Ex{\AA.11. {\bf (Double point  geometry)}}  Let $\f = dx_1\wedge\cdots\wedge dx_n +
 dy_1\wedge\cdots\wedge dy_n$  in $\bbr^{2n}$ for $n\geq 3$.  The only $\f$-planes are those parallel to the $x$ or $y$ axes.  An easy calculation shows that $dd^\f f=0$ 
 if and only if $f(x,y)=g(x)+h(y)$ for harmonic functions $g$ and $h$. However,  a function 
 $f(x,y)$ is $\f$-pluriharmonic if and only if it is harmonic in $x$ and $y$ separately.
 This is a simple example where (\AA.4) fails.

\medskip
In all  the following examples \fh functions are quite scarce.

\Ex{\AA.12. {\bf (Special Lagrangian geometry)}}   Consider the Special Lagrangian calibration
$\f={\rm Re}( dz)$ on $\bbc^n$.  For this calibration one can show that (\AA.4) is valid.  
Consequently, Lei Fu [Fu] has described 
all the $\f$-pluriharmonic functions.

\Prop{\AA.13} {\sl Let $f$ be a Special Lagrangian pluriharmonic function defined locally on $\bbc^n$, $n\geq 3$.  Then $f=A+Q$ where $A$ is affine and $Q$ is 
a traceless hermitian quadratic function.}

\pf  If $dd^\f f=0$ and $n\geq 3$ (so that $Z_{ij}$ and ${\overline Z}_{ij}$ are of different bi-degrees), then (\AA.1) implies that 
${\partial^2 f \over \partial {\bar z}_i \partial {\bar z}_j} = 0$ for all $i,j$.  Therefore, all third partial derivatives of $f$ are zero.  For polynomials of degree $\leq 2$ the result is transparent from (\AA.1).\qed

\Ex{\AA.14. {\bf (Associative, Coassociative and Cayley geometry)}}  
Consider one of the calibrations:\smallskip

1. (Associative) \qquad $\f(x\wedge y\wedge z) \ =\ \langle x, yz\rangle$ 
\qquad for $x,y,z \in {\rm Im} \bbo$

2. (Coassociative) \qquad\ $\psi \ =\ *\phi$  \qquad \qquad\qquad\qquad \ \ \ on ${\rm Im} \bbo$

3. (Cayley) \qquad $\Phi(x\wedge y\wedge z\wedge  w) \ =\ 
\langle x, y\times z\times w\rangle$ \qquad for $x,y,z, w \in  \bbo$

\smallskip\noindent
where $\bbo$ denotes the octonions.  As in the Special Lagrangian case
one can show that (\AA.4) is valid for each of these calibrations.   Furthermore, an application
of representation theory shows the maps $\BM_\f, \BM_{\psi}$ and $\BM_\Phi$ are injective.
 These calculations carry over to manifolds with  $G_2$ or   Spin$_7$-holonomy  to establish the following. 

\Prop{\AA.15}  {\sl Let $X$ be a manifold with holonomy contained in  $G_2$ or Spin$_7$
and  having dimension 7 or 8 respectively.  Suppose $\f$ is a parallel calibration
on $X$  of one of the three types above.  Then every  $\f$-pluriharmonic function
on $X$ is affine.  Moreover,  if the holonomy is exactly $G_2$ or Spin$_7$,
 every $\f$-pluriharmonic function is constant.  }
 
 \pf
 The first assertion follows because (\AA.4) is valid and the $\BM$-maps are injective.
 The second follows because any non-constant affine function on $X$ would reduce
 its  holonomy to a subgroup of $\{1\}\times SO_{n-1}$. \qed

\Ex{\AA.16. {\bf (Lie  group geometry)}}  Let $G$ be a compact simple Lie group with Lie algebra
$\gerG$, defined as the set of left-invariant vector fields on $G$.  
\smallskip
1)\  Consider the fundamental 
3-form $\f$ on $G$ defined by $\langle x, [y,z]\rangle$ and normalized to have comass one.
Calculations indicate that in all but a finite number of cases non-constant pluriharmonic functions
do not exist, however there are $\f$-submanifolds, namely the ``minimal'' SU$_2$-subgroups
(cf. [B], [T], [Th]).
\smallskip

2) Consider $*\f$.  The $*\f$-submanifolds are given by certain components of the cut locus $C$.
Is $G -C$ strictly $\f$-convex?

\Ex{\AA.17. {\bf (Gromov manifolds)}}    By a {\sl Gromov manifold} we mean an ensemble
$(X,\o, J, \langle\cdot,\cdot\rangle)$  where $(X,\o)$ is a symplectic manifold, $J$ is an 
almost complex structure on $X$ and $\langle\cdot,\cdot\rangle$ is a riemannian metric
with the property that
$$
\o(v,w)\ =\ \langle Jv, w \rangle
$$
for all $v,w\in T_xX$ at all $x\in X$. Every symplectic manifold has many Gromov structures.
Generically the almost complex structure $J$ is not integrable, and the only $\o$-pluriharmonic
functions are the constants. However, there are generally many $\o$-submanifolds (the pseudo-holomorphic curves) and there are many $\o$-plurisubharmonic functions as we shall see below.
It is important to note here that the operator $dd^c$ is not appropriate for this context since
$\nabla\o\neq0$.  
However, our notion of plurisubharmonicity works well and has the property that 
$\o$-plurisubharmonic functions are subharmonic on all pseudo-holomorphic curves.
 
 We note that on a Gromov manifold there exists a class of {\sl Lagrangian plurisubharmonic functions}
 with many good properties.
 For example, they are subharmonic  when restricted to any Lagrangian submanifold which is minimal.
 This is explored in a separate paper [HL$_5$].

\vfill\eject



\centerline{\bf \ZZ. \ The $\f$-Hessian.}\medskip

In this section we prove Theorem \AA.4 and Proposition \AA.7.  The arguments will involve ideas and notation
important for the rest of the paper.  A generalization of Theorem 1.4  to submanifolds which are {\sl $\f$-critical} can be found in Appendix A.

Recall (cf. [ON], p. 86) that the {\bf Hessian}, or second covariant derivative,   of a smooth function $f$ on a riemannian
manifold $X$ is defined on tangent vector fields
$V,W$ by
$$
\Hess(f)(V,W) \ \equiv\ V( W  f) - (\nabla_V W) f
\eqno{(\ZZ.1)}
$$
where $\nabla$ denotes the riemannian connection.
Note that $V( W  f) - (\nabla_V  W) f = V(\langle  W, \nabla f\rangle) - \langle \nabla_V  W, \nabla f\rangle = \langle W, \nabla_V(\nabla f)\rangle$ so that at a point $x\in X$, the Hessian is the symmetric 
2-tensor, or the symmetric linear map of $T_xX$ given by
$$
\Hess(f)(V) \ =\   \nabla_V(\nabla f).
\eqno{(\ZZ.2)}
$$
In terms of local coordinate vector fields 
$$
 \Hess(f) \left({\partial\over \partial x_i},  {\partial\over \partial x_k}\right) \ =\ 
 {\partial^2 f \over \partial x_i\partial x_j}  - \sum_k \Gamma_{ij}^k{\partial f\over \partial x_k}
$$
where $ \Gamma_{ij}^k$ are the standard Christoffel symbols of the riemannian connection.

Let $V$ be a real inner product space.
Given an element $\f\in \L^p V^*$, we define a linear map, central to this paper, 
$$
{\BM}_{\f} : \End(V)\ \arr\  \L^p V^*
\eqno{(\ZZ.3)}$$
by
$$
{\BM}_{\f} (A) \ \equiv \ {\Der}_{A^t}(\f)
$$
where ${\Der}_{A^t}$ denotes the extension of the transpose $A^t:V^*\to V^*$ 
to ${\Der}_{A^t}:\L^pV^*\to \L^pV^*$ as a derivation. That is, on simple vectors, one has
$$
{\Der}_{A^t}(v_1\wedge \dots \wedge v_p)\ =\  \sum_{k=1}^p v_1\wedge\dots\wedge A^t ( v_k ) \wedge 
\dots \wedge v_n
$$
\Note{} Recall that the natural inner product on $\End(V)$ is given by:
$$
\langle A,B   \rangle \ =\ \tr A B^*  \qquad {\rm for }\  A,B \in \End(V)
$$
Using this inner product we have  the adjoint map
$$
{\BM}_{\f}^* : \L^p V^*\ \arr\   \End(V)
\eqno{(\ZZ.4)}
$$
which will also be important.

\Note{}  If we identify $\End(V)$ with the Lie algebra ${{\hbox {\fr gl}}}(V)$ of $GL(V)$, then $\l_\f$ is the 
differential of the standard representation of $GL(V)$ on $\L^pV^*$ at $\f$. Therefore,  $\ker(\l_\f)$
is the Lie algebra of the subgroup $H_\f \equiv\{g\in GL(V) : g(\f)=\f\}$ and 
 $\ker(\l_\f) \cap {\rm SkewEnd}(V)$ is the Lie algebra of
 the compact subgroup  $K_\f  = H_\f \cap  {\rm O}(V)$.

\Def{\ZZ.1}  The {\bf $\f$-Hessian} of a function $f\in C^\infty(X)$ is the $p$-form $\CH^\f(f)$
defined by letting  the symmetric endomorphism    $\Hess f$ act on $\f$ as a derivation, i.e., 
$$
\CH^\f(f)   \ \equiv \    \Der_{\Hess^t f}( \f).
\eqno{(\ZZ.5)}
$$
In terms of the bundle map ${\BM}_{\f}  : \End(TX)\to \L^pT^*X$,
$$
\CH^\f(f)   \ \equiv \  {\BM}_{\f} (\Hess f)
\eqno{(\ZZ.6)}
$$
is the image of the Hessian of $f$.
\medskip

The second order differential operators $dd^\f$ and $\CH^\f$ differ by a pure 
first order operator. This is the first of the three equations needed to prove
Theorem \AA.4.

\Theorem{\ZZ.2}  {\sl  If $\f$ is a closed form on $X$, then}
$$
\CH^\f(f) \ =\ dd^\f f- \nabla_{\nabla f} \f
\eqno{(\ZZ.7)}
$$
\pf
By (\ZZ.2)  we have  $(\Hess f)(V) = \nabla_V \nabla f = [V, \nabla f] +\nabla_{\nabla f} V$, i.e.,
$$
\Hess f\ =\   -\cl_{\nabla f} + \nabla_{\nabla f}
$$
as operators on vector fields ($\cl$ = the Lie derivative).
The right  hand side of this formula has a standard extension to all tensor fields
as a derivation that commutes with contractions.  It is zero on functions, that is, it is
a bundle endomorphism whose value on $T^*X$ is minus the transpose of its value
on $TX$.  In particular, we find that
  ${\Der}_{\Hess^t f} = \cl_{\nabla f} - \nabla_{\nabla f}$ on $p$-forms, i.e., 
$$
\CH^\f(f) \ =\ \cl_{\nabla f}(\f) - \nabla_{\nabla f}\f
\eqno{(\ZZ.8)}
$$
Finally, since $d\f=0$, the classical formula    $d\circ \hk + \hk\circ d = \cl$   gives
$$
\qquad\qquad\qquad\qquad\qquad\qquad
dd^\f f \ =\ d(\nabla f\hk \f)\ =\ \cl_{\nabla f}(\f)  \qquad\qquad\qquad\qquad\qquad
\vrule width5pt height5pt depth0pt
$$
\medskip

Many of the nice results for the $dd^\f$-operator continue to hold in the 
non-parallel case after replacing it with the $\f$-Hessian. Perhaps even more importantly,
many properties of the $dd^\f$-operator  in the parallel case can best be understood by
considering the $\f$-Hessian. 

The second formula needed for the proof of Theorem \AA.4 is algebraic in nature,
involving  the bundle map ${\BM}_{\f}  : \End(TX)\to \L^pT^*X$.  Consequently, as before,
we replace $T_xX$ by a general inner product space $V$. If $\x$ is a $p$-plane
in $V$ (not necessarily oriented),  let $P_{\x}:V \to \x$ denote orthogonal projection.
The following, along with its reinterpretations (\ZZ.9$)'$ and (\ZZ.12), is a central result of this paper.

\Theorem{\ZZ.3}  {\sl  Suppose $\f$ has comass one.  For each $A\in \End(V)$,
$$
({\BM}_{\f}  A)(\x) \ =\ \langle  A, P_\x \rangle \qquad {\rm if\ } \x \in G(\f).
\eqno{(\ZZ.9)}
$$
Equivalently,}
$$
({\BM}_{\f}^*)(\x) \ =\ P_\x  \qquad {\rm if\ } \x \in G(\f).
\eqno{(\ZZ.10)}
$$
\medskip
Note that if $e_1,...,e_p$ is an orthonormal basis for the $p$-plane $\x$, then
$$
\langle  A, P_\x \rangle \ =\ \sum_{j=1}^p  \langle  e_j, Ae_j \rangle.
$$
Consequently, it is natural to refer to  $\langle  A, P_\x \rangle $  as the {\bf $\x$-trace of $A$}
and to use the notation 
$$
\tr_\x A\ \equiv\  \langle  A, P_\x \rangle.
$$
In particular, for each $A\in \End(V)$,
$$
({\BM}_{\f}  A)(\x) \ =\ \tr_\x A  \qquad {\rm if\ } \x \in G(\f).
\eqno{(\ZZ.9)'}
$$

Suppose $\x \in G(p,V)\subset \L_p V$ is a unit simple $p$-vector.  If
$a,b$ are unit vectors in $V$ with $a\in \span \x$ and $b\perp
\span\x$, then  
$$
b\wedge(a\hk\x)
$$
is called a {\bf  first cousin} of $\x$.  The first cousins of $\x$  
span the tangent space to the Grassmannian 
$G(p,V)\subset \L_p V$ at the point $\x$. 
Since $\f$ restricted to $G(p,V)$ is a maximum on $G(\f)$, this fact implies the
following  result, which we shall use frequently.

\Lemma{\ZZ.4. {\bf (The First Cousin Principle)}} {\sl  If $\phi\in    \L^p V^*$ has comass one and
$\x \in G(\f)$, then 
$$
\phi(\eta)\ =\ 0
$$
for all first cousins $\eta = b\wedge(a\hk\x)$ of $\xi$.}

\medskip

Note that $\Der_{(b\otimes a)^t} \f =\Der_{a\otimes b} \f = a\wedge(b\hk \f)$ and $\Der_{b\otimes a} \x =  b\wedge(a\hk\x)$
so that if $A=b\otimes a$ is rank one, then
$$
{\BM}_{\f}  (b\otimes a)(\x)\ =\ (\Der_{a\otimes b} \f )(\x)\  = \ \f(\Der_{b\otimes a} \x )
\ =\ \f(b\wedge(a\hk\x))
\eqno{(\ZZ.11)}
$$

\noindent
{\bf Proof of Theorem \ZZ.3.}  Pick an orthonormal basis for $\x$ and extend to an orthonormal basis of $V$.  It suffices to prove (\ZZ.9) when   $A=b\otimes a$
with $a$ and $b$ elements of this basis. It is easy to see that  $\langle  b\otimes a, P_\x \rangle =0$ 
unless $a=b\in\x$, in which case $\langle  a\otimes a, P_\x \rangle =1$.
By equation (\ZZ.11) we have ${\BM}_{\f}  (b\otimes a)(\x) = \f(b\wedge(a\hk\x))$ and
$b\wedge(a\hk\x)=0$ unless $a\in \x$ and either $b\in\x^{\perp}$ or $b=a$.
If $b\in\x^\perp$, then $(b\wedge(a\hk\x)$ is a first cousin of $\x$ and $\f((b\wedge(a\hk\x))=0$ by
the First Cousin Principle.  If $a=b\in\x$, then $b\wedge(a\hk\x)=\x$ and therefore
$\f((b\wedge(a\hk\x))=\f(\x)=1$. \qed

\medskip

Theorem \ZZ.3 has many   consequences.  We mention several. From (\ZZ.9$)' $ we have:

\Cor{\ZZ.5} {\sl  Suppose $(X,\f)$ is a calibrated manifold.  For each function $f\in C^\infty(X)$,}
$$
{\CH}^\f(f)(\x)\ =\ \tr_{\x}(\Hess f)  \qquad {\rm if\ }  \x\in G(\f).
\eqno{(\ZZ.12)}
$$

\medskip
This equation (\ZZ.12) is the second equation needed in the proof of Theorem \AA.4.

\Remark{} Equation  (\ZZ.12) provides an alternative definition of $\f$-plurisubharmonic
(as well as strictly $\f$-plurisubharmonic and $\f$-pluriharmonic) functions, which bypasses
the bundle map ${\BM}_{\f} $ and uses only the trace of the Hessian of $f$ on $\f$-planes $\x$.
\medskip

Another application of Theorem \ZZ.3 is given by:

\Cor{\ZZ.6} {\sl   If $A\in\End(V)$ is skew, then the $p$-form ${\BM}_{\f}  A$ vanishes on $G(\f)$.}
\medskip

See Remarks A.5 and A.7 for an extension of this to a recent result in [R].

Theorem \ZZ.3 has another useful consequence used  to prove  Lemma \AA.5.  
Note that for $A,B\in \Sym^2(V)
\subset\End(V)$, if   $A\geq 0$, $B \geq 0$, then $\langle A,B\rangle =\tr AB  \geq0$.
Hence for all $\x\in G_p(V)$ one has  $\langle e \otimes e, P_\x\rangle \geq0$, and more generally
 $\langle A, P_\x\rangle \geq0$ whenever $A\geq0$.  Since $df$ and $\nabla f$ are metrically equivalent, 
 $$
 {\BM}_{\f} (\nabla f \otimes \nabla f) = df\wedge ( \nabla f \hk \f) = df \wedge d^\f f.
\eqno{(\ZZ.13)} $$
 Therefore,  Theorem \ZZ.3 has the following consequence.  

\Cor{\ZZ.7} {\sl  For any $f\in C^\infty(X)$, }
$$
(df\w d^\f f)(\x)  \ =\ |\nabla f \hk \x |^2          \  \geq \ 0 \qquad{\sl for\ all\ }  \ \x \in G(\f).
\eqno(\ZZ.14)
$$

Theorem \ZZ.3 can also be used to understand the relationship between convex functions and \fp functions.  A function $f\in C^\infty(X)$ is called {\bf convex} if $\Hess f\geq 0$ at each point,
and it is called {\bf affine} if $\Hess f \equiv 0$ on $X$.  (If $f$ is affine, $\nabla f$ splits $X$ locally as a riemannian product
$\bbr\times X_0$.)

\Cor{\ZZ.8}  {\sl  Every convex function is 
$\f$-plurisubharmonic, and every strictly convex function is strictly \fp  (and every affine function is $\f$-pluriharmonic).}

\medskip

\Remark{\ZZ.9} The converse always fails; there are always \fp functions which are not convex.
To see this, consider first  the euclidean case with $X=V$ and $\f$ parallel.  
Recall that the orthogonal projections $P_e$ onto lines in $V$ generate  the extreme
rays of the convex cone of convex functions (positive semi-definite quadratic forms) in $\Sym^2V\subset \End(V)$.  This cone is self-dual. The projections $P_\x={\BM}_{\f}^*(\x) $ for $\x\in G(\f)$
generate a  proper convex subcone  (in fact a proper convex subcone of the cone generated by orthogonal projections onto $p$-planes).  Hence, by the Bipolar Theorem  there must exist a non-convex quadratic function $Q\in \Sym^2 V$ with $\langle Q, P_\x\rangle \geq 0$ for all $\x \in G(\f)$.  By (\ZZ.9), $Q$ is \ $\f$-plurisubharmonic.       
(Recall that for a convex cone $C \subset \bbr^n$ with vertex at the origin, the Bipolar
Theorem states that $(C^0)^0 = \overline C$ where $C^0 \equiv \{w \in \bbr^n : \langle w,v\rangle \geq 0$
for all $v\in C\}$ is the dual cone.)

\medskip

We now recall some elementary facts about submanifolds.
Given a submanifold $\overline X\subset X$, let $(\bullet)^T$ and 
$(\bullet)^N$ denote orthogonal projection of $T_xX$ onto the tangent and
normal spaces of $\ol X$ respectively. Then the canonical riemannian
connection $\overline \nabla$ of the induced metric on $\overline X$ is
given by ${\ol \nabla }_{ V}  W = (\nabla_V W)^T$ for tangent vector fields
$V,W$ on $\ol X$. The {\bf  second fundamental form} is defined by
$$
B_{V,W} \ \equiv \ (\nabla_V W)^N \ =\ \nabla_V W -\ol\nabla_V W.
$$
This is a symmetric bilinear form on $T\ol X$ with values in the normal
space. Its trace $H={\rm trace}\, B$ is the {\bf  mean curvature vector field}
of $\ol X$, and $\ol X$ is called a {\bf  minimal submanifold } if $H\equiv0$. 
Finally, let $\overline{\Delta}$ denote the Laplace-Beltrami operator   on $\ol X$
and  $\ol{\Hess}$ denote the Hessian operator on  $\ol X$.  The proof of the following 
is straightforward.
 $$
\ol \Hess(f)(V,W) = \Hess(f)(V,W) - B_{V,W} \cdot f
$$ 
for tangent vectors $V,W\in T \ol X$.  Taking the $T\ol X$-trace yields:
$$
\ol{\D} f \ =\ \tr_{T\ol X} \Hess f  - H(f).
$$

With a change of notation, this is the   final formula needed to prove Theorem \AA.4.

\Prop{\ZZ.10}  {\sl  Suppose $M$ is a $p$-dimensional submanifold of $X$ with mean curvature vector field
$H$.  Then for each $f\in C^\infty(X)$, }
$$
{\D_M} f \ =\ \tr_{TM} \Hess f  - H(f)  \qquad  {\ \rm on\ } M.
\eqno{(\ZZ.15)}
$$

\Cor{\ZZ.11} {\sl Suppose $M$ is a $\f$-submanifold of $X$. Then}
$$
\ch^\f(f)\bigr|_M\ =\ (\D_Mf) {\rm vol}_M.
\eqno{(\ZZ.16)}
 $$
\pf Combine (\ZZ.12) and (\ZZ.15) with the fact that $H=0$.\qed\medskip

Combining this with   (\ZZ.7) gives equation (\AA.2) and proves
Theorem \AA.4.

Given vectors $u,v\in V$, define $u\circ v \in \Sym^2(V)$ by
$u\circ v(w)=\half( \langle v,w \rangle  u + \langle u,w \rangle  v)$.

 \medskip
\noindent
{\bf Proof of Proposition \AA.7.}   
 In all cases $2\leq p\leq n$ 
$$
\Hess_x E\ =\ { p \over |x|^p  }\left( {1\over p}\cdot I - e\circ e  \right) \qquad {\rm with\ \ }e = {   x \over |x|   }.
$$
Set  $H  = {1\over p}I-e\circ e$.  Then  
$$
\l_\f(H) \ =\ \f - e\wedge (e \hk \f)\ =\ e\hk(e\wedge \f).
$$
Since $(e\wedge (e \hk \f))(\x) = \f (e\wedge (e \hk \x))$ and $e\wedge (e \hk \x)$ is a simple 
$p$-vector of norm $\leq 1$, for each unit simple $p$-vector $\x\in G(\phi)$  we have
$$
\l_\f(H)(\x) \ \geq \ 0  \qquad{\rm for\ all\ \ }  \x \in G(p,\bbr^n).
$$
This proves that E is \fp for all calibrations $\f$.

Finally, suppose $\x \in G(\f)$, i.e.,  $\f(\x)=1$.
Then  we have  $\l_\f(H)  =0$ if and only if $(e\wedge (e \hk \f))(\x) =1$,
which is equivalent to $(e\wedge (e \hk \x)) = \x$ or $e\in\span \x$.
This proves that $E$ is \pfp on $\bbr^n-\{0\}$ if and only if every vector $e$ is contained in a $\f$-plane.\qed

\medskip
For future reference we add a remark.

\Remark{\ZZ.12} When $\phi$ is harmonic,  the operator $d^\f$ can be expressed in terms of the Hodge $d^*$-operator as
$$
d^\f f \ =\  - d^*(f\f)
$$
and therefore
$$
dd^\f f\ =\  - dd^*(f\f).
$$
To prove this, first note that if $v\in T_xX$ and $\a\in T^*_xX$ are metrically equivalent, then
$v\hk \f =  (-1)^{(n-p)(p-1)} *(\a\wedge *\f)$.  Hence, $d^\f f = \nabla f \hk \f = 
 (-1)^{(n-p)(p-1)}*(df \wedge *\f) =   (-1)^{(n-p)(p-1)}\{  * (d(f*\f)-f(d*\f)) \}
=    (-1)^{(n-p)(p-1)} \{*d*(f\f)-f*d*\f \}$, and since $d^* = (-1)^{np+n+1} *d*$, we conclude that
$$
d^\f f\ =\  f d^*\f - d^*(f\f)
$$
so that the first equation holds if $\f$ is a harmonic form, and in particular if $\f$ is parallel.  Note also that for $\psi=*\f$
$$
d^\psi f\ =\ \pm *d(f\f) \and dd^\psi f \ =\ \pm *d^*d(f\f).
$$

\Remark {\ZZ.13. (Examples)} Corollary \ZZ.5 gives us the following basic fact.
\medskip
\centerline {  \sl  $f$ is $\f$-plurisubharmonic \ \ $\iff \ \ \tr_\xi \Hess\, f \ \geq\ 0 \quad \forall \, \x\in G(\f)$.}
\medskip\noindent
This infinitesimal  version of Theorem \AA.4 gives insight into the condition
of $\phi$-plurisubharmonicity. Consider for example the calibration 
$$
\Omega_p = \smfrac1{p!}\o^p
$$
where $(X,\o)$ is a K\"ahler, or more generally a Gromov, manifold.  In this case 
$$
G(\Omega_p) \ =\ G_{\bbc}(p, TX)\  \subset \  G_\bbr(2p, TX) 
$$
is exactly the set  of {\bf complex $p$-planes} in $TX$.
Thus a function $f$ is  $\Omega_p$-plurisubharmonic if and only if    $\tr_\xi  \Hess\, f \geq0$ for all
complex $p$-planes $\xi$.  In the K\"ahler case this means that  {\sl  $f$ is  $\Omega_p$-plurisubharmonic
if and only if it is subharmonic on all $p$-dimensional complex submanifolds.}
This condition can be expressed somewhat differently as follows. Any symmetric endomorphism
$A:TX\to TX$ can be decomposed as $A=A_{\rm sym}+ A_{\rm sk}$ where
$A_{\rm sym}= \half (A-JAJ)$ is hermitian symmetric with   real eigenvalues $\l_1\leq \l_2\leq \cdots \leq \l_n$.
If $\xi$ is a complex $p$-plane, then $\tr_\xi A = \tr_\xi A_{\rm sym}$.
Moreover, the infimum of $\tr_\xi A_{\rm sym}$ over such planes is $2(\l_1+\cdots+\l_p)$.
From this it follows that 
a function {\sl $f$ is  $\Omega_p$-plurisubharmonic if and only if  the eigenvalues $\l_1\leq  \cdots \leq \l_n$
of its hermitian symmetric hessian $\{\Hess\,f\}_{\rm sym}$ satisfy}
$$
\l_1+\cdots+\l_p \ \geq\ 0 \qquad {\sl at\ all\ points\ of\ } X.
$$
There is a parallel story on a quaternionic K\"ahler manifold $(X,\Psi)$ where
$\Psi$ is locally of the form $\Psi = {1\over 6}\{\o_I^2+\o_J^2+\o_K^2\}$ for
 orthogonal almost complex structures $I, J, K$ satisfying the standard relations.
For the calibration
$$
\Psi_p \ \equiv\ \smfrac{1}{(2p+1)!}\{\o_I^2+\o_J^2+\o_K^2\}^p, \quad  {\rm one\ has\ that}
$$
$$
G(\Psi_p) \ =\ G_{\bbh}(p, TX)\  \subset \  G_\bbr(4p, TX) 
$$
is the set of {\bf quaternionic $p$-planes}  in $TX$. Each symmetric endomorphism $A$ 
has a quaternionic hermitian symmetric part $A_{\rm qsym} = {1\over 4}(A-IAI-JAJ-KAK)$
with real eigenvalues  $\l_1\leq  \cdots \leq \l_n$, and a function $f$ is $\Psi_p$-\psh\ if and only if
the eigenvalues of $\{\Hess\, f\}_{\rm qsym}$ satisfy $\l_1+\cdots+\l_p \ \geq\ 0$ at each point.

In Section \FF \ we show that on a general $(X,\phi)$ the squared distance to any $\phi$-free submanifold $M$ is
strictly $\phi$-plurisubharmonic on a neighborhood of $M$.  This constructs
huge families of $\phi$-plurisubharmonic functions with topologically interesting
level sets.

\vskip .3in

\centerline{\bf  \KK. Elliptic Calibrations}\medskip

  In this paper we primarily restrict attention to $C^\infty$-functions.
However,  in this section we give the foundations for a theory of more general \fp functions 
developed in [HL$_4$].  
Suppose $(X,\f)$ is a calibrated manifold and assume that throughout this section that
$G(\f)$ is a fibre bundle over $X$.

\Def{\KK.1}  Given an everywhere  positive definite section $A$ of $\Sym^2(TX)$, the associated 
differential operator
$$
\D_A f\ \equiv \ \langle\Hess f, A\rangle
\eqno{(\KK.1)}
$$
will be called a {\bf Laplacian} on $X$.\medskip

The standard riemannian Laplacian is associated with the identity section of 
$\Sym^2(TX)$.  For each Laplacian $\D$ on $X$, classical potential theory is applicable.

\Def{\KK.2}   A calibration $\f$ is said to be {\bf mollified by a Laplacian $\D_A$}  
(or {\bf subordinate to a Laplacian $\D_A$}) if
every $f\in \fpsh$ is $\D_A$-subharmonic, i.e., $\D_Af\geq0$.
\medskip

\Ex{1}  Suppose that  $\f$ is a parallel calibration on $\rn$ and that
$\f=\sum_{j=1}^N \a_j\x_j$  is a positive linear combination  of  $\x_j\in G(\f)$. 
Assume that $\Sym^2(\rn)$ has only one irreducible component of dimension 1 
(the span of the identity $I$) under the subgroup of O$(n)$ that fixes $\f$.
Then $\f$ is mollified by   $ \D_A$ with $A=\BM_{\f}^* \f$.

\pf
If $f$ is \fp then
$$
\D_A f \ =\ \langle \Hess f, \BM_{\f}^* \f \rangle 
\ =\ \left \langle \Hess f, \sum_{j=1}^N\a_j \BM_{\f}^* \x_j \right\rangle
 \ =\   \sum_{j=1}^N\a_j \ch^\f(f)(\x_j)\ \geq\ 0.
$$
By the hypothesis $\BM_{\f}^*\BM_{\f}(I) =cI$ for some constant $c$.  For any calibration $\f$, one has
$\BM_{\f}(I)=p\f$.  Hence, $A=\BM_{\f}^*\f = {c\over p}I$.  Finally, $p^2\langle\f,\f \rangle 
= \langle\BM_{\f}(I), \BM_{\f}(I)\rangle = \langle  \BM_{\f}^*\BM_{\f}(I), I \rangle = cn$, proving that
$A=\BM_{\f}^* \f = {p\over n}|\f|^2\cdot I$  is positive definite.

Moreover, $\D_A = {p\over n}|\f|^2\D$ where $\D$ is the standard Laplacian.

\Remark{}  In general, the operator $\langle \Hess f, \BM_{\f}^* \f \rangle$ is not useful, as $\f$ may not be a positive combination of elements $\x_j\in G(\f)$.

\Def{\KK.3}     A calibration $\f$ is said to be {\bf elliptic} and  $G(\f)$ is said to {\bf involve all the variables}
if for every tangent vector $v\neq0$, there exists a $\f$-plane $\x$ with $v\hk \x\neq 0$.

\Ex{2} Let $\f=dx_1\wedge dy_1+\lambda dx_2\wedge dy_2$, with $0< \l  < 1$ on $\bbr^4$.
Then $G(\f) =  \{\x_0\}$, a single point in $G(2,\bbr^4)$, so that $\f$ is {\sl not} elliptic.
Note that the 2-form $\f$ {\sl does} involve all the variables.

\Ex{3} Suppose    $\f=dx_1\wedge\cdots \wedge dx_n+  dy_1\wedge\cdots \wedge dy_n
\in\L^n \bbr^{2n}$ with $n\geq 3$.  Then $G(\f)$ consists of only two points in $G(n,\bbr^{2n})$,
 but $\f$ is elliptic.  A function $f(x,y)$ is \fp if and only if it is separately subharmonic in $x$ and $y$.
 An example is $f(x,y)=u(x) v(y)$ with $u,v\geq0$ and subharmonic.

\Ex{4}  Suppose    $\Psi= \smfrac13(\half\o^2_I+\half\o^2_J+\half\o^2_K)$ on $\bbh^n$ is the 
quaternion calibration (cf. Example \AA.9).  Then the standard Laplacian $\D$ is a mollifying Laplacian
since
$$
\D f\ =\ \sum_{j=1}^n \tr_{\x_j} \Hess f\ =\ \sum_{j=1}^n(\ch^\Psi f)(\x_j)
$$
where $\x_1,...,\x_n$ are the axis $\bbh$-lines.  However, $\Psi$ is not a positive linear combination of elements in $G(\f)$.

\Theorem{\KK.4}  {\sl There exists a mollifying Laplacian for $\f$ if and only if $\f$ is elliptic.}

\Theorem{\KK.5}  {\sl Suppose $\f$ is elliptic.  A function $f\in C^\infty(X)$ is \fp 
if and only if $f$ is $\D$-subharmonic for every  mollifying Laplacian  $\D$. }
\medskip
The proofs of these two results will be given in the case $\f\in \L^pV$ is a parallel calibration on euclidean $n$-space.
Arguments for the more general case are essentially the same.

Let $\cp_+$ denote the convex cone in $\Sym^2(V)$ on $\{P_\x : \x\in G(\f)\}$.  That is, $A\in \cp_+$ 
if and only if 
$$
A\ =\ \sum_{j=1}^N \l_j P_{\x_j} \qquad {\rm with\ \ } \l_j>0 \ {\rm and\ } \x_j\in G(\f).
\eqno{(\KK.2)}
$$
Let $\cp^+$ denote the polar cone.  That is, $H\in \cp^+$ if and only if
$$
\langle H, P_\x\rangle\ =\ \tr_\x H\ \geq\ 0 \ \ {\rm for\ all\ }\x\in G(\f)
\eqno{(\KK.3)}$$
The Bipolar Theorem (cf. [S]) states that $A\in \Sym^2(V)$  can be expressed as in 
(\KK.2) if and only if 
$$
\langle H,A\rangle\ \geq\ 0 \ \ {\rm for\ all\ }  H \in \cp^+
\eqno(\KK.4)
$$

\Lemma{\KK.6}  {\sl Given $A\in \Sym^2(V)$, the associated operator $\D_A$ is a mollifying Laplacian for $\f$ if and only if 
\smallskip

1) \ \ $A\ =\ \sum_{j=1}^N \l_j P_{\x_j} \qquad {\rm with\ \ } \l_j>0 \ {\rm and\ } \x_j\in G(\f)$, and
\smallskip

2)\ \ $\x_1,...,\x_N$ involve all the variables, i.e., $v\hk \x_j=0, j=1,...,N$ implies $v=0$.
}
\pf Suppose 1) and 2) are valid.  Then $\langle Av,v\rangle = \sum \l_j |v\hk \x_j|^2$ and 2) implies that $A$ is positive definite.  Moreover, 
$$
\D_A f\ =\  \langle \Hess f,A\rangle \ =\ \sum \l_j \tr_{\x_j} \Hess f
\eqno{(\KK.5)}
$$
so that  if $f$ is \fp, then $f$ is $\D_A$-subharmonic.  

Conversely, suppose
$\D_A$ is a mollifying Laplacian for $\f$.  Take $f$ to be a quadratic function with $H=\Hess f
\in\cp^+$, so that $f$ is $\f$-plurisubharmonic.  Then $\D_A f=\langle H, A\rangle \geq 0$ for all such $H$.
As noted above (\KK.4) implies  (\KK.2).  
Finally, note that $\langle Av,v\rangle = \sum_{j=1}^N \l_j |v\hk \x_j|^2$
and therefore, since $A$ is positive definite 2) is verified. \qed

\medskip
\noindent
{\bf Proof of Theorem \KK.4}  Suppose there exists a mollifying Laplacian $\D_A$ for $\f$.  Then by 
Lemma \KK.6 1), we have $A=\sum_{j=1}^N \l_j P_{\x_j}$, and by 2) we have that 
given  $v\neq0$, there exists $\x_j\in G(\f)$ with $v\hk \x_j\neq0$.  Thus $\f$ is elliptic.

Conversely, if $\f$ is elliptic, then by compactness, there exists a finite number of $\x_j\in G(\f)$
such that
$$
A\ =\ \sum_{j=1}^N P_{\x_j} \in \cp_+
$$
is positive definite, thereby verifying 1) and 2). \qed

\medskip
\noindent
{\bf Proof of Theorem \KK.5} Suppose $f\in C^\infty(X)$  is subharmonic for every mollifying Laplacian
$\D_A$. Suppose $\x\in G(\f)$ and $A=\sum_{j=1}^N \l_j P_{\x_j}$ is a mollifying Laplacian.
Then  $A(t) = tA+(1-t)P_\x, \ \ 0<t\leq 1$ also determines a mollifying Laplacian by Lemma \KK.6.  Hence,
$$
\D_{A(t)} f \ =\ t \sum_{j=1}^N \l_j \tr_{\x_j} \Hess f + (1-t) \tr_{\x}\Hess f\ \geq\ 0.
$$
Taking the limit as $t\to0$, we obtain 
$
\tr_{\x}\Hess f\ \geq \ 0.
$
\qed

\bigskip
\centerline{\bf Generalized \fp Functions}
\medskip

Throughout this subsection we assume that $\f$ is an elliptic calibration.  The differential operator
$\Hess f$ extends to distributions $f$ on $X$ via duality producing a well-defined 
distributional section $\Hess f$ of the bundle $\Sym^2(TX)$.
By definition, a distributional section of a vector bundle $E\to X$ is a continuous linear functional on the space of smooth compactly supported sections of
$E^*\otimes \L^nT^*X$, or equivalently, on the space of $\wt s \equiv s\otimes *1$ for 
$s\in \G_{\rm cpt} E^*$.

\Def{\KK.7}  A distribution $f$ on $X$ is \fp if 
$$
(\tr_\x \Hess f)(\l)\ \geq\ 0
$$
for every smooth section $\x$ of $G(\f)$ and every smooth compactly supported non-negative multiple
$\l$ of the volume form on $X$.

\Ex{}  The fundamental function $E(x)$ in Proposition \AA.7 defines an
 $\lloc(\rn)$ function and, hence, is a distribution on $\rn$. It is \fp for any calibration $\f$ on $\rn$ of degree $p$.  This is easy to prove since (for $p<n$) the distributional hessian
 $$
 \Hess E \ =\ {p\over |x|^p}\left( {1\over p}I-{x\over|x|}\circ  {x\over|x|}   \right) \qquad{\rm on\ \ } \rn
 $$
has $\x$-trace
$$
\tr_\x\Hess E\ =\  {p\over |x|^p}\left( 1-{|x\hk \x|^2\over|x|^2}   \right).
$$

\medskip
Theorem \KK.5 extends from functions $f\in C^\infty(X)$ to distributions $f\in \cd'(X)$.   

\Theorem{\KK.8}  {\sl  Suppose $\D_A$ is a mollifying Laplacian for $\f$.  If $f$ is a \fp distribution,
then $\D_A  f \geq0$ is a non-negative measure, i.e., $f$ is $\D_A$-subharmonic.
Conversely, if $\D_A  f \geq0$ for each mollifying Laplacian $\D_A$, then $f$ is a \fp distribution.}
\pf
Assume that $A$ is of the form $A= \sum_{j=1}^N \l_j P_{\x_j}$ with $\l_j>0$ smooth and each $\x_j$ a smooth section of $G(\f)$. Then
$$
\D_A f\ =\ \langle A, \Hess f\rangle \ =\ \sum_{j=1}^N \l_j
\tr_{\x_j} \Hess f
$$
is a well defined distribution on $X$, and, by hypothesis, it pairs with every smooth, compactly
supported  non-negative multiple
$\l$ of the volume form, to give $(\D_A f)(\l)\geq 0$.  Hence $\D_A f $ is a non-negative regular Borel measure on $X$.
The proof of the converse is similar to that of Theorem \KK.5. \qed

\medskip

This theorem has a multitude of corollaries, deducible from the classical potential theory
for $\D_A$.  We list just two of the facts.
\smallskip

1)  Each $\D_A$-subharmonic distribution (and therefore each \fp distribution) belongs to $L^1_{\rm loc}(X)$, the space of locally Lebesgue integrable functions on $X$.
\smallskip

2) Each $\D_A$-subharmonic distribution (and therefore each \fp distribution) has a canonical
{\bf  classical} representative defined by
$$
f(x)\ =\ \lim_{r\to0}{1\over |B_r(x)|} \int_{B_r(x)} \, f\, d{\rm vol}
$$
which is $[-\infty,\infty)$-valued and upper semi-continuous on $X$. Here $B_r(x)$ denotes the 
ball of radius $r$ about $x$ and $ |B_r(x)|$ denotes its volume.

See [HL$_4$] for a development of upper semi-continuous \fp functions using these results.

 \vfill\eject


\centerline{\bf  \BB. Convexity in Calibrated Geometries}\medskip

We suppose throughout this section that $(X,\f)$  is a non-compact, connected  calibrated manifold
and all \fp functions are of class $C^\infty$.
\medskip

\Def{\BB.1}  If $K$ is a compact subset of $X$, we define the {\bf $(X,\f)$-convex hull of } $K$ by
$$
\wh K\ \equiv\ \{x\in X : f(x) \leq \sup_K f \ \ {\rm for\ all }\ f \in  \PSH(X,\f)\}.
$$
If $\wh K = K$, then $K$ is called {\bf $(X,\f)$-convex.}

\Lemma{\BB.2}  {\sl Suppose $K$ is a compact subset of $X$.  Then $x\notin \wh K$ if and only if there exists a smooth non-negative  \fp function $f$ on $X$ which is 
identically zero on a neighborhood of $K$ and has $f(x)>0$.  
Furthermore, if there exists a \fp function on $X$ which is strict at $x$, then  $f$ can be chosen to be strict at $x$.}

\pf  Suppose $x\notin \wh K$. Then there exists $g\in \PSH(X,\f)$ with $\sup_K g <0< g(x)$.
Pick $\varphi\in C^\infty(\bbr)$ with $\varphi\equiv 0$ on $(-\infty, 0]$ and with $\varphi>0$ and convex increasing on $(0,\infty)$.  Then $f=\varphi\circ g$ satisfies the required conditions
(See Lemma \AA.5a).
Furthermore, assume $h\in \PSH(X,\f)$ is strict at $x$. Then take  $\ol g = g +\epsilon h$.
For small enough $\epsilon$,  $\sup_K \ol g <0<\ol g(x)$.  If $\varphi$ is also strictly increasing on $(0,\infty)$, then $f=\varphi \circ \ol g$ is strict at $x$.  \qed

\medskip
\noindent
{\bf Note:} One sees easily that $\wh{\wh K} = \wh K$.  Therefore, if $\wh K$ is compact,  the function $f$ in Lemma \BB.2 can be taken
to be zero on a neighborhood of $\wh K$ (since one can replace $K$ with $\wh K$).

\Theorem{\BB.3}  {\sl  The following two conditions are equivalent.
\smallskip
1)  \ \ If $K\subset\subset X$, then $\wh K \subset\subset X$.

\smallskip

2) \ \  There exists a \fp  proper exhaustion function $f$ on $X$.  }

\Def{\BB.4} If the equivalent conditions of Theorem \BB.3 are satisfied, then $(X,\f)$ is a 
{\bf convex calibrated manifold } and {\bf $X$ is \fc.}

\smallskip
\noindent
{\bf Proof that  $2) \Rightarrow 1)$.}
If $K$ is compact, then $c\equiv \sup_K f$ is finite and $\wh K$ is contained in the compact pre-level set $\{x\in X : f(x) \leq c\}$.

\medskip
\noindent
{\bf Proof that  $1) \Rightarrow 2)$.}  A \fp proper exhaustion function on $X$ is constructed as follows.
Choose an exhaustion of $X$ by compact  $(X,\f)$-convex subsets 
$K_1\subset K_2\subset K_3\subset \cdots $ with $K_m\subset K^0_{m+1}$ for all $m$.
By Lemma \BB.2 and the compactness of $K_{m+2}-K_{m+1}^0$, there exists a \fp function $f_m\geq0$ on $X$ with  $f_m$ identically zero on a neighborhood of $K_m$   and $f_m>0$ on $K_{m+2}-K_{m+1}^0$.
By re-scaling we may assume $f_m>m$ on $K_{m+2}-K_{m+1}^0$. The locally finite sum 
$f=\sum_{m=1}^\infty f_m$ satisfies 2).\qed

\Lemma{\BB.5} {\sl Condition 2) in Theorem \BB.3 is equivalent to the {\rm a priori} weaker condition:
\smallskip

$2)'$  There exists a continuous proper exhaustion function $f$ on $X$ which is smooth and 

\quad \   \fp outside a compact subset of $X$.

\smallskip
\noindent
In fact if  $f$ satisfies $2)'$, then $f$ can be modified on a compact subset to be \fp on all of $X$.
Consequently, if $f$ satisfies $2)'$  and is strict outside a compact set, then its modification is also strict outside a compact set.}

\pf    For large enough $c$, $f$ is smooth and  \fp outside the compact set    $\{x\in X : f(x)\leq c-1\}$.  Pick a convex increasing function $\varphi\in C^\infty(\bbr)$   with $\varphi \equiv c$
 on a neighborhood of $(-\infty, c-1]$ and $\varphi(t)=t$ on $(c+1, \infty)$.  Then by Lemma \AA.5
the composition   $\varphi\circ f$ is \fp on all of $X$ (in particular smooth) and equal to $f$ outside 
  of the compact set  $\{x\in X : f(x)\leq c+1\}$.\qed

\Theorem{\BB.6} {\sl The following two conditions are equivalent:\smallskip

1)\ \ $K\subset\subset X  \ \Rightarrow \ 
\wh K\subset\subset X $, and $X$ carries a strictly \fp function.
\smallskip

2)\ \ There exists a strictly \fp proper exhaustion function for $X$.}

\Def{\BB.7} If the equivalent conditions of Theorem \BB.6 are satisfied, then   $(X,\f)$ is a
{\bf strictly convex calibrated manifold} or $X$ is {\bf strictly $\f$-convex}.
\medskip

\noindent
{\bf Proof of Theorem \BB.6.}  Suppose that $X$ is equipped with both a \fp proper exhaustion function
 $f$ and a strictly \fp function $g$. Then the sum $f+e^g$ is a strictly \fp exhaustion function.
 Now Theorem \BB.6 follows  immediately from Theorem \BB.3.\qed

\medskip

We shall construct many $\f$-convex manifolds in the course of our discussion
(See, in particular, \S \FF). However, we present some elementary examples here.

\Ex{1} Suppose $\f\in\L^p\bbr^n$ is a parallel calibration on $\bbr^n$.  Let
$f(x) = {1\over2}\|x\|^2$.  Then $dd^\f f=p\f$ and hence $f$ is a strictly \fp exhaustion.
That is,  $(\bbr^n, \f)$ is a strictly convex calibrated manifold.

\Ex{2}    Suppose $\f= dx_1\wedge\cdots\wedge dx_n$ on a domain $X$ in $\bbr^n$. 
Then $dd^\f f=(\D f)\f$ and $f$ is \fp if and only if $f$ is subharmonic.  Recall that if 
$K\subset\subset X$, then $\wh K = K \cup \{{\rm all\ the\  ``holes"\ in\ } K $ relative to $X\}$,
(connected components of $X-K$ which are relatively compact in $X$).
Thus $(X,\f)$  is strictly convex for any open set $X\subset \bbr^n$.

\medskip

It is instructive to extend this elementary example.

\Ex{3}      Suppose $\f= dx_1\wedge\cdots\wedge dx_p$ on a domain $X$ in $\bbr^n$
with coordinates $(x_1,...,x_p, y_1,...,y_{n-p})$. 
A function $f\in C^\infty(X)$ is  \fp if and only if $\D_x f\geq 0$ on $X$.
For a set $K\subset \bbr^n$, let $K_y$ denote the horizontal slice $\{x\in \bbr^p : (x,y)\in K\}$
of $K$.  Suppose that for each $y\in \bbr^{n-p}$, the horizontal slice $X_y$ has no
holes in $\bbr^p$.  Then $(X,\f)$ is strictly convex.  To prove this fact, it suffices to exhaust
$X$ by compact sets $K$ with the same property and show that each such $K$ 
is equal to its $(X,\f)$-hull.  Suppose $z_0=(x_0,y_0) \in X-K$.  Since $x_0$ is not in a hole
of $K_{y_0}$ in $\bbr^p$, we may choose (by Example 2) an entire subharmonic function
$g(x)$ with $g(x_0) >>0$ and $\sup_{K_{y_0}} g <<0$.  Now pick $\psi\in C^\infty_{\rm cpt}(\{y: |y-y_0|<\epsilon\})$ with $0\leq \psi\leq 1$ and $\psi(y_0)=1$.  Then $f(x,y) = g(x)\psi(y)$ is \fp 
and $f(z_0)=g(x_0)>>0$.  For $\epsilon$ sufficiently small, $\sup_K f\leq0$.  This proves  
$z_0$ does not belong to the $(X,\f)$-hull of $K$.

\Ex{4}  Let $\f=dx$ in $\bbr^2$ and set $X=\{(x,y) : x^2-c<y< x^2, \ |x|<1\}$.  Then $X$ is not $\f$-convex.  The closure of the hull of the compact subset
$K = ([-\e,\e]\times \{-\e\})\cup(\{\pm\e\}\times [-\e,0])$ of $X$ is easily seen to contain the origin.
Similarly, a domain of ``U''-shape, whose upper boundary along the bottom has a flat segment, is not $\f$-convex even though it is locally $\f$-convex
 (by Example 3).
\medskip

It is important to ``weaken''   this notion of strict convexity.

\Theorem{\BB.8}  {\sl  The following two conditions are equivalent:\smallskip

 1)\ \ $K \subset \subset X  \ \Rightarrow \  \wh K\subset\subset X $, and there exists a strictly \fp function defined outside a compact subset of $X$  \smallskip

2)\ \ There exists a   \fp proper exhaustion function  on $X$ which is strict outside a compact subset
of $X$.}

\Def{\BB.9}  If the equivalent conditions of Theorem \BB.8 are satisfied,  then the calibrated manifold  $(X,\f)$ is  
{\bf strictly convex at  $\infty$} or $X$ is {\bf strictly $\f$-convex at $\infty$.}
 
 \Remark{}  This is not  the standard terminology used in complex geometry where such spaces are called ``strongly (pseudo) convex".

\medskip
\noindent
{\bf Proof of Theorem \BB.8.}  Obviously 2) implies 1).  We will prove that 1) implies the following weakening of 2).
\smallskip

$2)'$ There exists a continuous proper exhaustion function $f$ on $X$ which is smooth and strictly \fp outside a compact subset of $X$.

\smallskip
By  Lemma \BB.5, Condition $2)'$ implies Condition 2).

Now assume 1).  Since $K\subset \subset X $ implies $\wh K\subset\subset X $, we know from Theorem \BB.3 that there exists a \fp exhaustion function $f$ for $X$.  Let $g$ denote the strictly \fp function which is only defined outside of a compact set.  We can assume this compact set is $\{x\in X : f(x)\leq c\}$ for some large $c$.  Then $h\equiv \max\{f+e^g, c\}$ is a continuous proper exhaustion function which, outside the compact set $\{x\in X : f(x)\leq c\}$, is strictly \fp (in fact, equal to$f+e^g$).
This proves $2)'$ and completes the proof of the theorem.\qed

\Cor{\BB.10} {\sl $(X,\f)$ is strictly convex at $\infty$ if and only if Condition $2)'$ holds.}

\vskip .3in
\centerline{\bf Cores.}\medskip

In each non-compact calibrated manifold $(X,\f)$ there are certain distinguished subsets
which play an important role in the $\f$-geometry of the space.  (In complex manifolds
which are strongly pseudoconvex, these sets correspond to the compact exceptional
subvarieties.)  The remainder of this section is devoted to a discussion of these subsets.

Given a function $f\in\fpsh$, consider the closed set
$$
W(f) \ \equiv\ \{x\in X: f {\rm \ is\ partially\ } \f-{\rm pluriharmonic\ at\ } x\}
$$
That is, $W(f)$ is the complement of the set
$$
S(f)\ \equiv\ \{x\in X : f \ \ {\rm is \  strictly \ }
\phi-{\rm plurisubharmonic  \ at\ \  } x \}
$$
Note that 
$$
W(\lambda f + \mu g) \ \subseteq\ W(f) \cap W(g) 
$$
for $f,g \in \fpsh$ and $\lambda, \mu >0$.

\Def{\BB.11}  The {\bf core} of $X$ is defined to be the intersection
$$
\Core(X) \ \equiv\ \bigcap W(f)
$$
over all $f\in\fpsh$. The {\bf inner core} of $X$ is defined to be the set ${\rm InnerCore}(X)$
of points $x$ for which there exists $y\neq x$ with the property that $f(x)=f(y)$
for all $f\in\fpsh$.

\Prop{\BB.12}  \qquad InnerCore$(C) \ \subset \ \Core(X)$.

\pf If $x\notin \Core(X)$, then there exists $g\in \PSH(X,\f)$ with $g$ strict at $x$.
Suppose $y\neq x$. Then if $\psi$ is compactly supported in a small neighborhood
of $x$ missing $y$, and $\psi$ has sufficiently small second derivatives, one has 
$f=g+\psi\in \PSH(X,\f)$. Obviously for such $f$, the values $f(x)$ and $f(y)$ can be
made to differ, so therefore $x\notin $ InnerCore$(X)$. \qed

\Prop{\BB.13}  {\sl  Every compact $\f$-submanifold $M$ without boundary in $X$ is contained in the inner core.}

\pf  Each $f\in\PSH(X,\f)$ is subharmonic on $M$ by Theorem \AA.4.  
Hence,   $f$ is constant on $M$. \qed

\Prop{\BB.14}  {\sl  Suppose $X$ is $\f$-convex.  Then $\Core(X)$ is compact
if and only if $X$ is strictly $\f$-convex at $\infty$,
and $\Core(X)=\emptyset$ if and only if $X$ is strictly \fc.}

\pf
If  $X$ is strictly $\f$-convex at $\infty$, then choosing $f$ to satisfy 2) in Theorem \BB.8, we see that the Core$(X) \subset W(f)$ is compact.  Obviously, strict $\f$-convexity implies that $\Core(X)=\emptyset$.

Conversely, if $\Core(X)$ is compact, then in the construction of the \fp exhaustion function in the proof of Theorem \BB.3 we may choose $K_1$ to be the $\phi$-convex hull of 
$
 {\Core(X)}.
 $
Then by the definition of $\Core(X)$ and Lemma \BB.2, each of the functions $f_m$ in that proof can be chosen to be strictly \fp 
on $K_{m+2}-K^0_{m+1}$. Hence the exhaustion $f=\sum_m f_m$ is strictly \fp outside a compact 
set containing the core.\qed
\medskip

A slight modification of this construction gives the following general result.

\Prop{\BB.15} {\sl  Suppose  $X$ is strictly $\f$-convex at $\infty$, and $K\subset X$ is
a compact, $\f$-convex subset containing the core of $X$.  Let $U$ be any
neighborhood of $K$. Then  there exists a proper \fp exhaustion function $f:X\to \bbr^+$
 which is strictly \fp  on $X-U$, and  identically zero on a neighborhood of $K$.
}
\pf    
Choose $K_1$
as  in the construction of the \fp exhaustion function given in the proof of Theorem \BB.3.
Let $K(\epsilon)$ denote the compact $\e$-neighborhood of $K_1$. Then
$$
K_1 \ =\ \bigcap_{\e>0} \wh{K(\e)}.
\eqno{(\BB.1)}
$$
If $x\in \bigcap_{\e>0} \wh{K(\e)}$, then for each $f\in \PSH(X,\f)$, we have $f(x) \leq \sup_{K(\e)} f$.
However,   $\inf_{\e} \sup_{K(\e)} f = \sup_K f$, and we conclude that $x\in \wh K$.  Thus we can choose $K_2 \equiv  \wh K(\e)$ in our construction of $f$, and for  small enough $\e$ we have $K_2\subset U$ as well as $K_1\subset K_2^0$. The proof is now completed as in the proof of Proposition \BB.14.  \qed
\medskip

Obviously, many questions concerning
$$
{\rm InnerCore}(X)\ \subseteq\ \Core(X)\ \subseteq\ \wh{\Core(X)}
$$
remain to be answered.

\vskip .3in



\centerline{\bf   Examples of Complete  Convex Manifolds and Cores}

\medskip

In \S \FF \ (Theorem \FF.4) we shall show that there are many strictly $\f$-convex domains in any calibrated manifold
$(X,\f)$.  They can have quite arbitrary topological type within the strictures imposed by Morse Theory
and $\f$-positivity of the Hessian.  However, it is  also interesting geometrically to ask for convex manifolds which are complete.
 
 In fact, there exist enormous families of  complete calibrated manifolds $(X,\phi)$ with $\nabla\phi=0$ 
 which are  strictly $\phi$-convex at infinity.  For example, any  $(X,\phi)$ which is asymptotically locally euclidean   (ALE) is such a creature. In this case the radial function on the asymptotic cone at infinity
 is strictly convex. It appears likely the corresponding assertion also holds  for quasi ALE spaces.
  For the general construction of such spaces with  SU(n), Sp(n), G$_2$, or Spin$_7$ holonomy, the reader is referred to the book of Joyce [J].
  
Certain  manifolds of this type have been quite explicitly constructed, and in these cases
 one can explicitly construct \fp exhaustion functions and identify the cores. We indicate how
 to do this below.
 
 We begin however with an observation in dimension 4. Every crepant resolution of singularities
 of $\bbc^2/\G$, for a finite subgroup $\Gamma\subset {\rm SU}(2)$, admits a  Ricci-flat ALE K\"ahler metric.  On each such manifold there exists an $S^2$-family of parallel calibrations
$$
 \cc\ =\  \{u\omega + v\vf+\ w \psi : u^2+v^2+w^2=1\}
 $$
where $\o$ is the given K\"ahler form, 
 $\varphi= {\rm Re}\{\Phi\}$  and $\psi  ={\rm Im}\{\Phi\}$ and $\Phi$ is a parallel section of the canonical bundle $\kappa_X$.   Let $E=\pi^{-1}(0)$ be the exceptional locus of the resolution. Then for any $\phi\in \cc$  we have
 $$
 \Core(X,\phi) \ =\ \cases {E & if $\phi=\o$ \cr \emptyset & otherwise.}
 $$
 This follows from the fact that each $\phi\in \cc$ is in fact the K\"ahler form for a complex structure
 on $X$ compatible with the given metric.  With this complex structure $X$ is pseudo-convex, and by the Stein Reduction Theorem  (cf.   [GR, p. 221])  we know its core is the union of its compact complex subvarieties.  For $\phi\neq \o$ there are no such subvarieties since by the Wirtinger inequality  
 (cf.  [L$_{1,2}$]),  applied to $\phi$, they would necessarily be homologically mass-minimizing, and by the same result applied to $\o$  any such subvariety is $\o$-complex (and therefore a  component of $E$).
 
 \Ex {1. (Calabi Spaces)}  Let $X\to \bbc^n/\bbz_n$ be a crepant resolution of $\bbc^n/\bbz_n$
 where the action on $\bbc^n$ is generated by scalar multiplication by $\tau = e^{2\pi i/n}$.
Following Calabi [C] we define the function $F:\bbc^n/\bbz_n\to \bbr$ by
$$
F(\rho) \ =\ \root n \of {\rho^n+1} + {1\over n}\sum_{k=0}^{n-1} \tau^k \log\left(\root n \of {\rho^n+1} -\tau^k\right)
$$
where  $\rho \equiv \|z\|^2$ (pushed down to $\bbc^n/\bbz_n$), and the log is defined by 
choosing $\arg \zeta \in (-\pi,\pi)$.
 We then define a Kahler metric on $\bbc^n/\bbz_n-\{0\}$ by setting
 $$
 \o \ =\ {1\over 4} d d^c F.
 $$
 Calabi shows that {\sl this metric is Ricci-flat  and (when pulled back)  extends to a Ricci flat metric on $X$.  The parallel form $\Phi = dz_1\wedge\cdots\wedge dz_n$ extends to a parallel section of } $\kappa_X$.   This metric is given explicitly on $\bbr^{2n}/\bbz_n$ by
 $$
ds^2\ =\  F'(\rho) |dx|^2 + \rho F''(\rho) dr \circ d^c  r
$$
 where $r=\|x\|$.  Define $G(\rho)$ by setting $G'(\rho)=F'(\rho)+\rho F''(\rho)$ and $G(0)=0$.
 Then direct calculation shows that
 $$
 dd^\f G \ =\ 2n \phi
 $$
 where $\phi = {\rm Re} \{\Phi\}$.  Hence, {\sl  $X$ is a complete,  strictly $\phi$-convex manifold.}
 
 \Ex{2. (Bryant-Salamon Spaces)}  Let $P$ denote the principal Spin$_3$-bundle of $S^3$ and 
 $$
 S\ \equiv\ P \times_{{\rm Sp}_1} \bbh
 $$
 the associated spinor bundle, where $\bbh $ denotes the quaternions.
Bryant and Salamon have explicitly constructed  a complete riemannian metric with
G$_2$-holonomy on the total space of $S$. (See [BS, page 838, Case ii].)  Let $\rho=|a|$
for $a\in \bbh$ (pushed-down to $S$) and let $Z\subset S$ denote the zero section.  Then a direct calculation shows that the function
$$
F(\rho) \ =\ (1+\rho)^{{5\over 6}} \ \ {\sl  is\ strictly\ } \varphi-{\sl plurisubharmonic\ on \ } S-Z
$$
 where $\varphi$ denotes the associative calibration on $S$.   Since $Z$ is an associative  submanifold
 we conclude that
 $$
 \Core(S)\ =\ Z.
 $$
 In an analogous fashion the authors construct a complete riemannian metric with 
 Spin$_7$-holonomy on the total space $\wt S$ of a spinor bundle over $S^4$.
 (See  [BS, page 847, Case ii].)  A similar calculation shows that there exists an exhaustion
 function which is strictly $\Phi$-plurisubharmonic on  $\wt S - \wt Z$ where $\Phi$ denotes the Cayley
 calibration  $\wt Z$ the zero-section of $\wt S$. Since $\wt Z$ is a Cayley submanifold, we conclude that $$
 \Core (\wt S)\ =\ \wt Z.
 $$

 \vfill\eject


\centerline{\bf \CC. Boundary Convexity.}
 
\medskip

Suppose $\O\subset\subset X$ is an open set with smooth boundary $\bo$, where $(X,\f)$ is a            non-compact calibrated manifold.  A $p$-plane $\x\in G(\f)$ at a point  $x\in \bo$ will be called  {\bf tangential} if $\span\x\subset T_x\bo$.

  \Def{\CC.1}  Suppose that $\rho$ is a {\bf defining function} for  $\bo$,
  that is, $\rho$ is a smooth function defined on a neighborhood 
  of $\overline{\Omega}$ with $\Omega =\{x: \rho(x)<0\}$ and  $\nabla \rho \neq 0$ on $\partial \Omega$.
  If
  $$
 \ch^{\f} \rho(\xi) \ \geq\ 0 \ \ \  {\rm for\ all\  tangential\ } \xi\in G_x(\phi), \ x\in \bo,
   \eqno{(\CC.1)} 
 $$
  then  $\partial \Omega$ is  called {\bf $\phi$-convex}.
  If the inequality in (\CC.1) is strict for all $\xi$, then 
  $\partial \Omega$ is  called {\bf strictly $\phi$-convex}.  If 
  $\ch^{\f}(\xi) =0$ for all $\xi$ as in (\CC.1), then    $\partial \Omega$ is {\bf $\phi$-flat}.
  
  Each of these conditions is a local condition on $\bo$.  In fact:

  \Lemma {\CC.2} {\sl  Each of the three conditions in Definition \CC.1  is independent of the choice of defining function $\rho$.
  In fact, if $\overline{\rho}=u\rho$ is another choice with $u>0$ on $\partial \Omega$, then on $\bo$ }
  $$
 \ch^{\f}(\overline{\rho})(\xi)\ =\ u \ch^{\f} (\rho)(\xi) \ \ \ {\rm for \ all \  tangential\ }\ \ \xi\in G(\phi)
    \eqno{(\CC.2)}
    $$
    
  \pf
Since $\rho=0$ on $\bo$ and $\nabla\rho\perp\bo$, we have  by (\ZZ.12) that 
 $$\eqalign{
  \ch^{\f}(\overline{\rho})(\xi)\ &=\  \tr_\x \left[\Hess(u\rho)\right]  \
  \ =\   \tr_\x \left[ u \, \Hess( \rho)  +2\nabla u\circ \nabla\rho + \rho\,\Hess u       \right]   \cr
&=\  u\, \tr_\x (\Hess \rho) +2\tr_\x (\nabla u\circ \nabla\rho) + \rho\tr_\x(\Hess u)   
\ =\   u \, \tr_\x (\Hess \rho) = u\,\ch^{\f}(\rho). \ \ \vrule width5pt height5pt depth0pt  
}
 $$
  
  \Cor{\CC.3} {\sl  Assume $\f\in \L^p\bbr^n$ is a calibration.  Suppose $\bo$ is (strictly) $\phi$-convex in $\bbr^n$, and locally near a point $p\in \bo$, 
  let   $\bo$ be graphed over its tangent space by a function $x_n=u(x')$ for linear coordinates
  $(x',x_n)$ on $\bbr^n$.  Then each nearby hypersurface: $x_n = u(x')+c$ is also (strictly) $\phi$-convex.
  }\medskip
    The next lemma and its corollary will be used to establish the main results of this section.
  
  \Lemma {\CC.4} {\sl Suppose $\rho$ is a smooth real-valued function 
  on a riemannian manifold, and $\psi:\bbr\to\bbr$ is smooth on the image of $\rho$.
  Then
    $$
\tr_\xi \Hess \psi(\rho)  \ =\ \psi'(\rho)\tr_\xi \Hess\rho +   \psi''(\rho)|\nabla\rho\hk \xi|^2 
  \eqno{(\CC.3)}
  $$
 for all oriented tangent $p$-planes $\xi$.}
 
 \pf
 We first calculate that
 $
 \Hess \psi(\rho)  = \psi'(\rho)  \Hess\rho +   \psi''(\rho)\nabla \rho \circ \nabla \rho
 $
and then note that $\tr_\xi( \nabla \rho \circ \nabla \rho )= |\nabla\rho\hk \xi|^2 $.\qed

\Cor{\CC.5}  {\sl  With $\delta = -\rho$ and $\rho <0$, one has }
$$
\tr_\xi \Hess(-\log \delta)\ =\ {1\over \delta}\tr_\xi \Hess \rho + {1\over \delta^2}|\nabla \rho\hk \x|^2
\eqno{(\CC.4)}
$$
\pf 
Take $\psi(t)=-\log(-t)$ for $t<0$, and note that $\psi'(t)=-1/t $ and  $\psi''(t)=1/t ^2$, so that $\psi'(\rho)=
 1/\delta$  and $\psi''(\rho)= 1/\delta^2$.\qed

   We now come to the main result of this section.

  \Theorem{\CC.6}   {\sl   Let $\Omega\subset\subset X$ be a compact domain with strictly \fc boundary.          Suppose $\delta = -\rho$ is an arbitrary ``distance function'' for $\bo$, i.e.,   $\rho$
  is an arbitrary  defining function for $\bo$.
  Then $-\log \,\delta$ is strictly \fp   outside a compact     subset of $\O$. Thus, in particular, the domain $\Omega$ is strictly $\phi$-convex at $\infty$.   }
  $ \over $
  
 \pf   
Applying (\ZZ.12) to  Corollary \CC.5 shows that at each point $x\in\O$
 near $\bo$, we have  
  $$
 \ch^{\f} (-\log\, \delta)(\x)\ =\ {1\over {\delta}} \ch^{\f} (\rho)(\x) + {1\over {\delta^2}} |\nabla\rho\hk \x|^2
  \eqno{(\CC.5)}
  $$ 
 for all $\x\in G(\f)$.  Note that at $x\in \bo$,  $ |\nabla\rho\hk \x|^2$ vanishes if and only if $\x$ is tangential to $\bo$.    For notational convenience we set 
 $$
 \cos^2\theta(\x) \ =\ {|\nabla\rho\hk \x|^2\over |\nabla \rho|^2} \ =\ \langle P_{\span \nabla \rho}, P_{\span \x} \rangle.
$$
Then
the inequality $|\cos\theta|<\epsilon$ defines a fundamental neighborhood system for $G(p,T\bo)\subset G(p,TX)$.  By restriction     $|\cos\theta|<\epsilon$ defines a fundamental neighborhood system for $G(\phi)\cap G(p,T\bo)\subset G(\phi)$.  The hypothesis of strict  $\phi$-convexity for $\bo$ implies that there exists $\overline \epsilon >0$ so that
  $( \ch^{\f}  \rho)(\x) \geq \overline\epsilon$ for all $\phi$-planes $\x$  at points of $\bo$ with $|\cos\theta|<\epsilon$
  for some $\epsilon>0$. (Note that if there are no $\phi$-planes tangent to $\bo$ at a point
  $x$, then there are no $\phi$-planes with $|\cos\theta|<\epsilon$ for sufficiently small
  $\epsilon$ in a neighborhood of $x$.)   Consequently, we have by equation (\CC.5)  that 
   $$
  \ch^{\f}  (-\log\delta)(\x) \ \geq \ {\overline \epsilon \over 2\delta}
  $$
    near $\bo$ for all  $\phi$-planes $\x$  with $|\cos\theta|<\epsilon$.
    
  Now choose $M>>0$ so that $\ch^\phi(\rho)(\x)\geq -M$ in a neighborhood of $\bo$ for all 
  $\x\in G(\phi)$.  Then, by (\CC.5) 
 $$
   \ch^\phi (-\log\, \delta)(\x) \ \geq\   -{M\over \d}+  {1\over \d^2}  |\nabla \rho\hk \x|^2.
  $$
  If $|\cos\theta|\geq \epsilon$, this 
   is positive in a neighborhood of $\bo$ in $\Omega$.  This proves that $-\log \d$ is strictly \fp near $\bo$.    By Corollary \BB.10 the domain $\O$ is strictly $\f$-convex at infinity.
  \qed  
 
 \medskip
 
 Although a general defining function for a strictly \fc boundary may not be $\f$-plurisubharmonic, for some applications the following is useful.
 
 \Prop{\CC.7} {\sl  Suppose $\O\subset\subset X$ has strictly \fc boundary $\bo$ with defining function     $\rho$.  Then, for $A$ sufficiently large, the function $\overline \rho \equiv \rho + A \rho^2$ is strictly \fc in a neighborhood of $\bo$  and also a defining function for $\bo$.}

 \pf 
 By Lemma \CC.4  and (\ZZ.12)
 $$
 \ch^\f(\rho)(\x)\ =\ (1+2A\rho) \ch^\f( \rho)(\x) +2 A |\nabla \rho\hk \x|^2 
  \qquad{\rm for\ all\ } \x\in G(\f).
 \eqno{(\CC.6)}
 $$
  As noted in the proof  of Theorem \CC.6,  strict boundary convexity
  implies the existence of  $\e, \overline\e>0$ so that, along $\bo$, $ \ch^\f( \rho)(\x)\geq \overline \e$ if  $\x\in G(\f)$ with 
 $|\cos\theta(\x)|<\e$.  Therefore $ \ch^\f( \overline\rho)(\x) \geq (1+2A\rho)\overline \e$  for all such $\x$.  Choose a lower bound $-M$ for $\ch^\f( \rho)(\x)$ over all $\x\in G(\f)$ in  a neighborhood of $\bo$.  Then by (\CC.6), 
 $\ch^\f( \rho)(\x)\geq -(1+2A\rho)M + 2|\nabla\rho|^2 A \e^2$ for $\x\in G(\f)$ with 
 $|\cos \theta(\x)|\geq \e$ at points of $\bo$.  For $A$ sufficiently large, the right hand side is $>0$ in some neighborhood of $\bo$.\qed
 \medskip

One might hope for a converse to Theorem \CC.6, e.g.,   if the domain $\O$ is \fc then the  boundary
is $\f$-convex.  However, the following elementary example  shows that this is false. 

\Ex{}  Let $\f\equiv dx\wedge dy$ in $\bbr^3$ as in Example 3 of section \BB.  Let $X$ denote the solid torus  obtained by rotating the disk $\{(y,z) : y^2+(z-R)^2 < r^2\}$ about the $y$-axis.  Since each slice $X_z$ has no holes in $\bbr^3$, the domain $X$ is $\f$-convex
(cf. Example 3 of \S \BB).  However, the boundary torus $\partial X$ is \fc  if and only if  $2r\leq R$.  This follows from  an elementary calculation which uses the obvious defining function and   Definition  \CC.1  (or by using Proposition \CC.13 below)

\Qu{\CC.8}  For which strictly convex calibrated manifolds is it true that \fc sub-domains have \fc boundaries?  More generally, when is the $\f$-convexity of a domain a local condition at the 
boundary?
  \medskip
  
  A weak partial converse to Theorem \CC.6 is given by the following.
  
  \Prop{\CC.9}  {\sl  Suppose the calibration is parallel, and set $\d=\dist(\bullet, \bo)$ in $\Omega$.
  If $-\log \d$ is   \fp near $\bo$, then $\bo$ is \fc.
  }

\Note{\CC.10}  Examples show that the strict convexity of $-\log\d$ near $\bo$ is stronger than 
$\f$-convexity for $\bo$.

\pf   Set $\rho =-\delta$ on $\overline\O$   near  $\bo$.
  Suppose that $\bo$ is not $\phi$-convex. Then there exist $x\in \bo$ and 
  $\x_x\in G_x(\phi)$ with $\span(\x)\subset  T_x(\bo)$ and $( \ch^\phi \rho)(\x_x)<0$.
  Let $\gamma$ denote the geodesic segment in $\Omega$ which emanates orthogonally from $\bo$ at $x$. Since $\delta$ is the distance function, $\gamma$ is an integral curve of $\nabla \delta$.  Let 
  $\x_y$, $y\in \gamma$ denote the parallel translation of $\x_x$ along $\gamma$.  Then $\x_y$ is a 
  $\phi$-plane with $\span(\x_y)\perp \nabla \rho$ for all $y$.  By formula (\CC.5), since 
  $\nabla\rho \hk \xi_y=0$, we have
  $$
   \ch^\phi (-\log\, \delta)(\x_y) \ =\ {1\over \delta}  \ch^\phi(\rho)(\x_y) \ <\ 0
  $$
  for all $y$ sufficiently close to $x$.  Hence, $-\log \delta$ is not \fp near $\bo$.  \qed
  
\bigskip

\centerline{\bf The Second Fundamental Form}
\medskip

  The $\phi$-convexity of a boundary can be equivalently defined in terms of its second fundamental form.  Note that if $M\subset X$ is a smooth hypersurface with  a chosen unit normal field $n$ we have a quadratic form $II$ defined on $TM$ by 
  $$
  II(V,W) \ =\ \langle  B_{V,W}, n \rangle
  $$  
  where $B$ denotes  the second fundamental form of $M$ discussed in \S \ZZ.
  For example, when $H=S^{n-1}(r)\subset \bbr^n$ is the euclidean sphere of radius $r$, oriented by the outward-pointing unit normal,  we find that $II(V,W) = -{1\over r} \langle V,W\rangle$.

For the sake of completeness we include a proof of the following standard fact.

\Lemma{\CC.11} {\sl Suppose $\rho$ is a defining function for $\O$ and let II denote the second
fundamental form of the   hypersurface $\bo$ oriented by the outward-pointing normal.
Then 
$$
\Hess\, \rho\,\bigl|_{T\bo} \ =\ -|\nabla \rho| \, II
$$
and therefore
$$
\tr_\xi \Hess \rho\ =\ -|\nabla \rho |\, \tr_\xi II
$$
for all $\xi\in G(p, T\partial \O)$.}
\pf  Suppose $e$ is a tangent field on $\bo$. Extend $e$ to a vector field tangent to the level sets
of $\rho$. By definition $II(e,e)=\langle \nabla_e e, n\rangle$ where 
$n = \nabla \rho/|\nabla \rho|$ is the outward normal.
Then $(\Hess\, \rho)(e,e) = e(e \rho) -(\nabla_e e)\rho = -(\nabla_e e)\rho = -\langle \nabla_e e, \nabla \rho\rangle = - |\nabla \rho| \langle \nabla_ee, n\rangle$.
\qed

  \Remark{}
    Recall that a defining function $\rho$ for $\Omega$ satisfies  $|\nabla \rho| \equiv 1$ in a neighborhood of $\bo$ if and only if $\rho$ is the signed distance to $\bo$ ($<0$ in $ \Omega$ and $>0$ outside of $\Omega$).   In fact any function $\rho$ with $|\nabla \rho| \equiv 1$  in a riemannian manifold is, up to an additive constant,   the distance function to (any) one of its level sets.
  In this case it is easy to see that
      $$
    \Hess \rho \ =\ \left( \matrix{0 & 0 \cr 0 & - II}       \right)
   \eqno{(\CC.7)} 
   $$
   where $II$ denotes the second fundamental form of the hypersurface $H=\{\rho =\rho(x)\}$
   with respect to the normal $n= \nabla \rho$ and the blocking in (\CC.7) is with respect to the 
   splitting $T_xX = \span (n_x) \oplus T_x H$.     For example let $\rho (x) =\|x\|\equiv r$ in $\bbr^n$. 
   Then direct calculation shows that $\Hess\, \rho
 = {1\over r}(I-\hat x \circ \hat x)$ where $\hat x = x/r$.

  \Cor{\CC.12}  {\sl   For all tangential $\x\in G(\f)$}
  $$
  ( \ch^\f \rho)(\x) \ =\ -|\nabla \rho| \tr_\x II.
  $$
  \pf Apply Theorem \ZZ.3.\qed
  \medskip
  
  As an immediate consequence we have
  
  \Prop{\CC.13}  {\sl  Let $\Omega\subset X$ be a domain with smooth boundary $\bo$ oriented by
  the outward-pointing normal. Then $\bo$ is \fc if and only if its second fundamental form satisfies
  $$
  \tr_\x II \ \leq \ 0 
  $$
  for all $\phi$-planes $\xi$ which are tangent to $\bo$.  This can be expressed more geometrically by saying that
  $$
   \tr\left\{ B \bigr|_\xi\right\} \ \ \ {\sl must\ be\  inward-pointing\  }
  $$
  for all tangential $\phi$-planes $\xi$.
  }

  \Remark{}
  If $\rho $ is the signed distance to $\bo$, then equation (\CC.7) together with Lemma \CC.11
  can be used to simplify (\CC.4).  An arbitrary $p$-plane $\xi$ at a point can be put in the canonical form $\xi=(\cos\theta n + \sin \theta e_1)\wedge e_2\wedge \cdots \wedge e_p$ with $n=\nabla \rho$
  and $n, e_1,...,e_p$ orthonormal.  Then $\eta=e_1\wedge \cdots \wedge e_p$ is the tangential
  projection of $\xi$.  Note that $\tr_\xi \Hess \rho = -\sin^2\theta \tr_\eta II$ and that 
  $|\nabla \rho\hk \xi|^2=\cos^2\theta$, so that (\CC.4) becomes
  $$
  \tr_\xi \Hess (-\log \delta) \ =\ -{1\over \delta}\sin^2\theta\, \tr_\eta II +{1\over \delta^2} \cos^2\theta.
  $$

We finish this section with a useful characterization of strictly convex domains.

\Theorem {\CC.14}  {\sl  Let $(X,\f)$ be a strictly \fc manifold and $\O\subset\subset X$ a 
domain with smooth boundary.  Then the following are equivalent.
\smallskip

(i)  \ \ $\bo$ is strictly \fc.
\smallskip

(ii) \ $\tr_\x\left\{ II_{\bo}\right\} <0$ for all tangential $\f$-planes  $\x$.
\smallskip

(iii)   There exists a defining function $\rho\in  C^\infty(\overline\O)$ for $\bo$ which is strictly
\fp

\qquad  on a neighborhood of $\overline\O$.   } 

\pf
It is clear that (iii) $\Rightarrow$ (i) $\Leftrightarrow$ (ii), so we need only prove that 
 (i) $\Rightarrow$ (iii).  Suppose $\bo$ is strictly $\f$-convex. By Proposition \CC.7
 we may assume that $\bo$ has a defining function $\rho_0$ 
 which is strictly \fp in a neighborhood of $\bo$.  By the Inverse Function Theorem there 
 is a neighborhood $U$ of $\bo$ and a diffeomorphism:
 \smallskip
 \centerline{$\bo\times [- 2\e, 2\e]\harr{\cong} {}U$ \quad
such that\ \  $\rho_0(x,t)=t$. 
}
\smallskip \noindent
 Let $\rho_1:X\to\bbr^+$ be a strictly \fp proper exhaustion function.
 Choose $\d$ with $0<\d<<\e$.
By replacing $\rho_1$  with  $a\rho_1-b$ for suitable $a,b>0$ we may assume that 
\smallskip
\centerline{$-\e -\d \ < \  \rho_1\ < \  -\e+\d<0$ \qquad on the neighborhood $U$ }
\smallskip
\noindent
Note that $\max\{\rho_0,\rho_1\}=\rho_1$
in the region where $t\equiv \rho_0< -\e-\d$, and that  $\max\{\rho_0,\rho_1\}=\rho_0$ where 
$t\equiv \rho_0 >  -\e+\d$ (in particular, in a neighborhood of $\bo$).

By Remark \AA.6  the function $\max\{\rho_0,\rho_1\}$
can be approximated by a smooth strictly \fp function $\rho$ on $U$ which agrees 
with $\max\{\rho_0,\rho_1\}$ outside the  compact subset of $U$  where $|\rho_0-\rho_1|\leq \d$.
We see that $\rho=\rho_1$ when $t<-\e-2\d$ and $\rho=\rho_0 \ (=t)$ when $t > -\e+2\d$. 
Therefore, $\rho=\rho_0$ in a neighborhood of $\bo$, and $\rho$ extends smoothly to  $\O$ by setting $\rho=\rho_1$ in $\O-U$. \qed

\vfill\eject



\vskip .3in

\centerline{\bf \FF.  \   $\f$-Free Submanifolds and Topology}

\medskip

Somewhat surprisingly,  for any calibration $\f$ there is a precise integer bound on the 
homotopy dimension  of any strictly $\f$-convex domain. This is the first result below. 
After establishing it, we  show that on the other hand, subject to this bound, there exist strictly convex domains of almost arbitrary topological type.

\Def{\FF.1}   The {\sl free dimension}, denoted   $\fdim$,  of a calibrated manifold $(X,\phi)$ is the maximum 
dimension of a  linear subspace in $TX$ which contains no $\f$-planes.
Such subspaces will be called {\sl $\f$-free}.

\Theorem {\FF.2}  {\sl   Suppose $(X,\f)$ is a strictly $\f$-convex manifold. Then $X$ has the homotopy type of a CW complex of dimension $\leq \fdim$.}

\pf Let $f:X\to \bbr^+$ be a strictly \fp proper exhaustion function. By perturbing we may assume that 
$f$ has non-degenerate critical points. The theorem follows if we show that each critical point
has Morse index $\leq \fdim$ (cf. [M]). If this fails, then there is a critical point $x$ at which
 $\Hess_x f$ has at least  $\fdim+1$ negative eigenvalues.  In particular, there exists a subspace
$W\subset T_xX$  of dimension $= \fdim+1$ with  $\Hess_x f\bigr|_W <0$.   However, by definition of $\fdim$,  $W$ must contain a $\f$-plane $\x\in G(\f)$, and since $f$ is strictly \fc, we must have $\tr_\x \Hess_x f >0$, a contradiction.
\qed

\medskip\noindent {\AAA E\BBB XAMPLES:}

\item{(a)} If $(X,\o)$ is a K\"ahler manifold of real dimension $2n$, then ${\rm fd}(\o) = n$.

More generally one has ${\rm fd}({1\over p!}\o^p) = n-p+1$.

\item{(b)}If $(X,\vf)$ is a Ricci-flat  K\"ahler manifold (Calabi-Yau manifold) of real dimension $2n$
with Special Lagrangian calibration $\vf$, then ${\rm fd}(\vf) = 2n-2$.

\item{(c)}If $(X,\Psi)$ is a quaternionic K\"ahler manifold or
 hyperK\"ahler manifold of real dimension $4n$ with the quaternionic calibration, 
 then ${\rm fd}(\Psi) = 3n$.
More generally for the calibration $\Psi_p \equiv {1\over (2p+1)!}  (\o_I^2 + \o_J^2 + \o_K^2 )^p$
 one has ${\rm fd}(\Psi_p) = 3(n-p+1)$.

\item{(d)}If $(X,\f)$ is a 7-manifold with  an associative calibration $\f$, then ${\rm fd}(\f) = 4$.

\item{(e)}If $(X,\psi)$ is a 7-manifold with  a coassociative calibration $\psi$, then ${\rm fd}(\psi) = 4$.

\item{(f)}If $(X,\Phi)$ is an 8-manifold with  a Cayley calibration $\Phi$, then ${\rm fd}(\Phi) = 4$.

\medskip\noindent {\AAA C\BBB OMMENTS:}

For (a), note  that every real subspace  of dimension $n+1$ in $\bbc^n$ contains a complex line
and is therefore not free.
The free subspaces of dimension $n$ are exactly the {\bf totally real} $n$-planes -- those for
which $JW\cap W=\{0\}$.

For the second statement, recall that the $\phi$-planes are exactly the $J$-invariant subspaces of dimension $2p$.  Now if $W\subset \bbc^n$ has codimension $\leq n-p$, then
$\dim_{\bbr} \{W\cap JW\}\geq 2p$, and so $W$ is not 
$\phi$-free.  However for a generic $W$ of codimension $n-p+1$, the maximal complex subspace
of $W$ satisfies 
$\dim_\bbr \{W\cap JW\}  =   2p-2$.
\smallskip
For (b), we first show that every real hyperplane $H\subset \bbc^n$ 
contains a Special Lagrangian $n$-plane
and is therefore not free.  Choose a unit vector $n\perp H$ and consider 
the orthogonal decomposition $\bbc^n = (\bbr n) \oplus (\bbr Jn) \oplus H_0$
where $H_0 = H\cap J(H)$ is the maximal complex subspace of $H$.
If $L_0\subset H_0$ is a Lagrangian subspace of $H_0$, then $L=(\bbr Jn) \oplus L_0$
is a Lagrangian subspace of $H$.   Rotating  $L_0$ in $H_0$ makes
$L$ Special Lagrangian as claimed.
We now observe  that for a  real subspace $W\subset \bbc^n$ of   dimension $2n-2$,
$$
W \ \ {\sl is\  } \vf\ {\sl free\ } \quad \iff \quad J(W^\perp)\not\subset W   
\iff \quad W\ \ {\sl is \ symplectic, \ i.e.,}\ \ \o^{n-1}\bigr|W \ \neq\ 0.
$$
For the first equivalence note  that if $J(W^\perp)\subset W$,
then the construction above gives a Special Lagrangian $L\subset W$.
Conversely, given $L\subset W$, $J(L) = L^\perp = (L^\perp\cap W) \oplus W^\perp$
and so $J(W^\perp)\subset W$.
For the second equivalence, note that  $J(W^\perp)\subset W$  implies that $J(W^\perp)$
lies in the null space of $\o\bigr|_W$. Conversely, if $v\in W$ lies in the null space 
of $\o\bigr|_W$, then $J(\span\{v, Jv\} )   \subset W^\perp$

\smallskip
For (c) suppose $V\subset \bbh^n$ has codimension $<n$. Then $V\cap I(V)\cap J(V)\cap K(V)$ is a non-trivial quaternionic subspace of $V$ and so $V$ cannot be free.
A subspace  $W$ of real codimension $n$  is free if and only if 
$W^\perp\cap IW^\perp = W^\perp\cap JW^\perp= W^\perp\cap KW^\perp = \{0\}$.
 
 For the second statement we use that fact that $G(\Psi_p)$ is exactly the set of 
 quaternionic linear subspaces of quaternionic dimension $p$ in $\bbh^n$.
 The argument then proceeds as in part two of (a).

\smallskip
For (d), suppose $V\subset {\rm Im}{\bf O}$ has dimension  5. 
 Let $x,y$ be an orthonormal basis of $V^\perp$.  Choose any unit vector $\e\in V$ perpendicular 
 to $z\equiv x y$.  Then $\span\{z, \e, z\e\}\subset V$ is associative, and so $V$ is not free.
(To see that $z\e \in V$, note that left multiplication by $z$ is an isometry which preserves
the quaternion subalgebra span$\{x,y,z\}$ and therefore also preserves its orthogonal complement.)
 We now claim that 
 
 \smallskip
 \centerline{\sl  a 4-plane $W\subset {\rm Im}{\bf O}$
  is free if and only if $W^\perp$ is not $\phi$-isotropic,
 i.e., $\phi\bigr|_{W^\perp}\ne 0$.} \smallskip
 
 \noindent
 To see this, suppose there exists an associative
 3-plane $V\subset W$. Then $V^\perp = \bbr \e \oplus W^\perp$ is coassociative
 (where $\e\in W$ is a unit vector perpendicular to $V$). Choose an orthonormal 
 basis $x',y',z'$ of $W^\perp$. Coassociativity is equivalent to the fact that the
 4-form $(*\phi)(\e, x',y',z') = \langle \e, [x', y', z']\rangle = \pm 1$ where $[\cdot,\cdot,\cdot]$ is the
 associator. We now recall the general equality $\phi(x', y', z')^2 + \| [x', y', z']\|^2 = 1$ (cf.
 [HL$_1$]) from which it follows that $\phi(x', y', z')=0$, that is, $W^\perp$ is $\phi$-isotropic.
 Conversely, supposing $W^\perp$ is $\phi$-isotropic, the equality shows that 
 $[x', y', z']$ is a unit vector and therefore $U\equiv \span\{x', y', z',[x', y', z']\}$ is coassociative.
 Hence, $U^\perp \subset W$ is associative.
 
\smallskip
For (e), suppose $V\subset {\rm Im}{\bf O}$ has dimension
5.  Let $x,y$ be an orthonormal basis of $V^\perp$. 
Then $U=\span \{x,y, x y \}$ is associative, and so $U^\perp\subset V$ is coassociative.
Hence, $V$ is not free.  Of course the free 4-planes are exactly those which are not coassociative.

\smallskip
For (f), suppose $V\subset  {\bf O}$ has dimension 5.
  Let $x,y,z$ be an orthonormal basis of $V^\perp$. Then $W=\span \{x,y,z,x\times y\times z\}$ is a Cayley plane, and so $W^\perp\subset V$ is also Cayley.  Hence, $V$ is not free. 
  The free 4-planes are exactly those which are not Cayley.
\medskip

We now show that within the homotopy restrictions imposed by Theorem \FF.2, 
the possible topologies
for strictly $\f$-convex manifolds are vast.

Let $(X,\f)$ be a calibrated manifold.  A $p$-plane $\x$ is said to be {\sl tangential} to a 
submanifold $M\subset X$ if $\span \x \subset T_xM$.

\Def{\FF.3}  A closed submanifold $M\subset X$ is {\bf  $\f$-free} if there are no $\f$-planes $\x\in G(\f)$ which are tangential to $M$.
If the restriction of the calibration $\f$ to $M$ vanishes,   $M$ is called  {\bf  $\f$-isotropic}.

\medskip

Note that  $\f$-isotropic submanifolds are  $\f$-free.
Each submanifold of dimension strictly less than the degree of $\f$ is $\f$-isotropic and hence
automatically $\f$-free.  Furthermore, in dimension $p$ the generic local submanifold is $\f$-free.
Depending on the geometry, this may continue through a range of dimensions greater than $p$.

\Theorem {\FF.4}  {\sl Suppose $M$ is a closed submanifold of $(X,\f)$ and let 
$\dis(x) \equiv {1\over2}{\rm dist}_M^2(x)$ denote half the square of the distance to $M$.  Then
$M$ is $\f$-free  if and only if 
the function $\dis$ is strictly \fp at each point in $M$ (and hence in a neighborhood of $M$).
}

\pf  We begin with the following.

\Lemma{\FF.5}  {\sl  Fix $x\in M$ and let $P_N:T_xX\to N$ denote orthogonal
 projection onto the normal plane   of $M$ at $x$.  Then for each $\x \in G(\f)$ one has}
$$
\{\BM_\f(\Hess_x \dis)\}(\x) \ =\ \langle P_N, P_\x\rangle
\eqno {(\FF.1)}
$$
\pf  By  Theorem \ZZ.3
$$
\{\BM_\f(\Hess_x f)\}(\x) \ =\ \langle \Hess_x f, P_\x\rangle
\eqno {(\FF.2)}
$$
for any function $f$.
The lemma then follows from the assertion that 
$$
\Hess_x \dis \ =\ P_N.
\eqno{(\FF.3)}
$$
To see this we first note that the Hessian of any function $f$  can be written
$$
\Hess f (V,W) \ =\ \langle V, \nabla_W( \nabla f)\rangle 
\eqno{(\FF.4)}
$$
for all  $V,W\in T_xX$.
It follows that if $\nabla f =0$ on the submanifold $M$, then $T_xM \subset {\rm Null}(\Hess_x f)$.
Thus, with respect to the decomposition $T_xX=T_xM \oplus N$ we have
$$
\Hess_x \dis \ =\ \left(\matrix { 0&0\cr 0&A\cr }\right)
$$
and it remains to show that $A$ is the identity. To see this, set $\delta(x)= \dist_M(x)$  
and note that $\nabla \delta = n$ is a smooth unit-length vector field near (but not on)
$M$ whose integral curves are geodesics emanating from $M$.   Hence, $$\nabla_n (\nabla \dis) =
\nabla_n (\nabla {1\over2}\delta^2) = \nabla_n (\delta n) = n + \delta \nabla_n n= n.$$
Taking limits along normal geodesics  down to $M$ gives the result.
\qed
\medskip

Theorem \FF.4 now follows from the fact that 
$$
\langle P_N, P_\x\rangle\ \geq\ 0 \ {\rm with \ equality \ iff }\ \span\x \subset N^\perp = T_xM.
$$
\qed

The following result gives us a powerful, very general method for constructing 
strictly \fc domains in $(X,\f)$.

\Theorem{\FF.6}  {\sl  Suppose $M$ is a $\f$-free submanifold of $(X,\f)$.
  Then there exists a fundamental neighborhood system $\cf(M)$ of $M$ consisting of strictly
  $\f$-convex domains. Moreover,
  
  \smallskip
  
  (a)\  $M$ is a deformation retract of each $U\in \cf(M)$.

  \smallskip
  
  (b)\  \ ${\PSH}(V,\f)$ is dense in  ${\PSH}(U,\f)$ if $U\subset V$ and $V,U\in\cf(M)$.

  \smallskip
  
  (c)\  Each  compact set $K\subset M$ is ${\PSH}(U,\f)$-convex for each $U\in \cf(M)$.
  }

\Remark{} The existence of $\f$-free   submanifolds gives the existence of strictly \fc domains 
with essentially every topological type permitted by Morse Theory (Theorem \FF.2).  Note in particular
that if $M\subset X$ is  $\f$-free, then {\bf every submanifold of $M$ is also $\f$-free}.

\pf   We construct tubular neighborhoods of $M$ as follows.  Let  $\epsilon \in C^\infty(M)$ 
be a smooth function which vanishes at infinity and has the property that  for each $x
\in M$ the ball $\{y\in X:{1\over2}\dist(y,x)^2\leq \epsilon(x)\}$ is compact and geodesically convex.
Assume also that $\e$ is sufficiently small so that the exponential map gives a diffeomorphism
$$
\exp : N_\e \ \arr\ U_\e
$$
from the open set $N_\e$  in the normal bundle $N$ defined by
${1\over2} \|n_x\|^2<\e(x)$ to the  neighborhood
$$
U_\e \ =\ \{x\in X : \dis(x) < \e(x)\}.
\eqno{(\FF.5)}
$$
of $M$ in $X$.
Each $U_\e$ admits a deformation retraction onto $M$.

By Theorem \FF.4 the function $\dis={1\over2}{\rm dist}_M^2(\cdot)$ is strictly \fp\ on a neighborhood of $M$, which we can assume to be $W$.  We impose the following additional condition on the function $\e\in C^\infty(W)$.
$$
\dis - t\e \ \ \  {\rm is\ strictly \ } \f-{\rm plurisubharmonic \ on\ } \ W \ \ {\rm for }\ \ 0\leq t \leq 1.
\eqno{(\FF.6)}
$$
Since (\FF.6) is valid as long as $\e$ and its first and second derivatives vanish sufficiently fast
at infinity, it is easy to see that the family $\cf(M)$ of neighborhoods $U_\e$ constructed above with $\e$ satisfying (\FF.6) is a fundamental neighborhood system for $M$.

Obviously, the function $\psi \equiv (\e-\dis)^{-1}$ is a proper exhaustion for $U_\e$.  Recall that if $g$ is a positive concave function, then $1/g$ is convex, or more directly, calculate that
$$
\Hess \psi \ =\ \psi^2\Hess(\dis-\e) +\psi^3\nabla(\e-\dis)\circ \nabla(\e-\dis).
\eqno{(\FF.7)}
$$
Applying $\BM_\f$ to (\FF.7) proves that $(\e-\dis)^{-1}$ is strictly \fp\ on $\{\dis<\e\}=U_\e$.
Hence, $U_\e$ is strictly $\f$-convex.

To prove parts (b) and (c) one uses Proposition \KK.16 in [HL$_2$], 
characterizing denseness of $\PSH(V,\f)$ in   $\PSH(U,\f)$ in terms of relative convexity, and argues exactly as on page 302 of 
[HW$_1$].       \qed

\medskip

\Ex{\FF.7}  As mentioned above, Theorem \FF.6 exhibits a rich family of \fc \ domains in $(X,\f)$.
For example, let $M\subset X$ be {\sl any submanifold of dimension} $< p = \deg \f$. Then by
\FF.6, $M$ has a fundamental system of neighborhoods each of which is a strictly $\f$-convex domain homotopy equivalent to $M$.

\Ex{\FF.8} Interesting examples occur in all the calibrated geometries examined in depth in [HL$_1$].
Suppose for instance that $X$ is a Calabi-Yau manifold with Special Lagrangian calibration
$\f$.  Then {\sl any complex submanifold $Y\subset X$ (of positive codimension)  is $\f$-free.}  It follows that any smooth submanifold of $Y$ is also $\f$-free.
\bigskip

\centerline{\bf Zero Sets  of Non-negative Strictly \fp\  Functions}\medskip

We now consider  the following two classes of subsets of $(X,\f)$.
\smallskip

\qquad (1)\ Closed subsets $A$ of $\f$-free submanifolds.
\smallskip

\qquad (2)\ Zero sets  of non-negative strictly \fp\ functions $f$.
\smallskip

These two classes are basically the same, as described in the following two propositions.

\Prop{\FF.9}  {\sl  Suppose $A$ is a closed subset of a $\f$-free submanifold $M$ of $X$.  Then
there exists a non-negative function $f\in C^\infty(X)$ with
\smallskip

(a) \ \ $A=\{x\in X: f(x)=0\}$
\smallskip

(b) \ \ $f$ is strictly \fp\ at each point in $M$ (and hence in a neighborhood 

\qquad  of $M$ in $X$).
}

\pf  Since $M$ is a closed submanifold, the function  $\dis$ in Theorem \FF.4 can be extended to
$h\in C^\infty(X)$ which agrees with $\dis$ in a neighborhood of $M$ and satisfies
$$
h\geq 0 \and \{h=0\}= M.
$$

Choose $\psi\in C^\infty(X)$ with $\psi\geq0$ and $A=\{x\in X: \psi=0\}$.  Now choose $\e\in
 C^\infty(X)$  with $\e(x)>0$ for all $x\in X$, and with $\e$ and its derivatives sufficiently small so that $f\equiv h+\e\psi$ is strictly \fp \ on $M$.
\qed

\Prop{\FF.10} {\sl  Suppose $f  \in C^\infty(X)$ is a non-negative function which is strictly 
\fp at each point in $A \equiv \{x\in X:f(x)=0\}$. Given a point $x\in A$ there exists a neighborhood $U$ of $x$ and a proper $\f$-free  submanifold $M$ of $U$ such that $A\cap U\subset M$.
}
\pf
Given $x\in A$ we may choose geodesic normal  coordinates $(z,y)$ in a neighborhood $U$ at $x$  so that
$$
\Hess_x f\ =\ \left( \matrix{0&0\cr0&\Lambda\cr}\right)\eqno{(\FF.8)}
$$
where $\Lambda$ is the diagonal matrix diag$\{\lambda_1,...,\lambda_r\}$, $r$ is the rank
of $\Hess_x f$, and $\lambda_j\neq 0$ for $j=1,...,r$.
Set 
$$
M\ =\  \left \{w\in U: {\partial f\over\partial y_1}=\cdots ={\partial f\over\partial y_r}=0\right\}.
$$
Since $\nabla{\partial f\over\partial y_1},..., \nabla{\partial f\over\partial y_r}$ are linearly independent at $x$, $M$ is a codimension $r$ submanifold locally near $x$.  

Note that ker$(\Hess_x f)=T_xM$. It remains to show that ker$(\Hess_x f)$ is $\f$-free
 (since  if $M$ is $\f$-free at $x$, then $M$ is $\f$-free in a neighborhood of $x$). This  is proved in Lemma \FF.11 below.\qed\medskip
 
 \Lemma{\FF.11}  {\sl  Suppose $f$ is  strictly \fp at $x\in X$.   Then $\ker(\Hess_x f)\subset T_xX$
 is $\f$-free.     }

\pf
If  $\ker(\Hess_x f)\subset T_xX$ is not $\f$-free, there exists $\x\in G(\f)$ with 
$(\Hess_x  f)\bigr|_{\span \x}=0$.  Consequently, 
$(\ch^\f f)(\x) = \BM_\f(\Hess_x  f)(\x) = \tr_\x( \Hess_x f)=0$, and $f$ is not strict at $x$.
\qed

\Remark{\FF.12}    Theorem \FF.6 can be generalized as follows.
Suppose $M=\{f=0\}$ is the zero set of a non-negative strictly \fp function $f$ on $(X,\f)$.
Then there exists a fundamental neighborhood system $\cf(M)$ of $M$ consisting of  strictly $\f$-convex domains which satisfy (c) of Theorem \FF.6.  The neighborhoods $U_\e\in\cf(M)$ are defined by 
$U_\e=\{x\in X: f(x)<\e(x)\}$ where $\e>0$ is a $C^\infty$ function on $X$ vanishing at infinity along with  its first and second derivatives so that $f-\e$ remains strictly \fp.  The proofs of (b), (c) and (d)
are essentially the same as in Theorem \FF.6.

\medskip

We conclude with  the following useful observation.

\Prop{\FF.13} {\sl Let $M$ be a submanifold of $(X,\f)$ and $f$ a smooth function defined on
a neighborhood of $M$ such that:
\smallskip

(1) $\nabla f \equiv 0$ on $M$, and \smallskip

(2)  $f$ is strictly \fp at all points of $M$. \smallskip

\noindent
Then $M$ is $\f$-free.}

\pf  By (\FF.4) we see that $TM\subseteq \ker(\Hess f)$ at all points of $M$.  We then apply Lemma \FF.11.\qed\medskip

\Cor{\FF.14} {\sl Let $f$ be a non-negative, real analytic function on $(X,\f)$ and consider the real analytic subvariety $Z\equiv \{f=0\}$.  If $f$ is strictly \fp at points of $Z$, then each stratum of $Z$
is $\f$-free.}

\vskip .3in

\centerline{\bf  Appendix A.  Submanifolds which are $\f$-critical.}\medskip

\def\s{A}

Here we establish a useful extension of Theorem \ZZ.3  to certain $\x$ which are not $\f$-planes.
Let $G \equiv  G(p,V)$ denote  the Grassmannian  of oriented p-planes in the inner product  space $V$, considered as the subset $G \subset \L_p V$ of unit simple vectors.

\Def{A.1}  Given $\f\in \L^pV^*$ an element $\x\in G$ is said to be a {\bf $\f$-critical
point} if $\x\in G$ is a critical point of the function $\f\bigl|_G$.
Equivalently, $\f$ must vanish on $T_\x  G \subset \L_pV$.
Let $G^{\rm cr}(\f)$ denote the set of $\f$-critical points.  \medskip

Note that if $\f$ is a calibration on $G$, i.e., $\sup \f\bigl|_G =1$, then
$$
G(\f)\ \subset \ G^{\rm cr}(\f)
$$
since for each $\x\in G(\f)$ the form $\f$ attains its maximum value 1 at $\x$.
Equation (\ZZ.9$)'$  extends from $G(\f)$ to $G^{\rm cr}(\f)$ as follows

\Prop{A.2}  {\sl   Suppose $\f\in \L^p V^*$  and   $\s \in\End(V)$.  Then for all $\x\in  G^{\rm cr}(\f)$}
$$
{\BM}_{\f} (\s)(\x)\ =\ (\tr_\x \s)\f(\x)
$$
\pf This is an immediate consequence of the more general Proposition A.4 below.\qed
\medskip

 Recall that at  a point $\x \in G$
there is  a canonical isomorphism:
$$
T_\x G\ \cong \ \Hom (\span \x, (\span \x)^\perp).
\eqno(A.1)
$$
On the other hand, $T_\x G$ is canonically a subspace of $\L_p V$.
It is exactly the subspace spanned by the first cousins of $\x$.
More specifically,  the isomorphism (A.1) associates to 
$L:\span \x \to (\span \x)^\perp$ the $p$-vector  $\Der_L \x$.

\Def{A.3} Let $\s\in\End (V)$ be a   linear map.  At each point $\x \in G$ we define
a tangent vector 
$$
\Der_{\wt A} \x\ \in\ T_\x G
$$
where  $\wt A =  P_{\x^{\perp}} \circ { \s }\circ P_\x$. This vector
field  $\x\to \Der_{\wt A} \x$  on $G$ is called the {\bf $A$-vector field}.
 
\Remark{} A straightforward calculation shows that if $A$ is symmetric, this  $A$-vector field on $G$      
is the gradient 
of the height  function $F_A : G \to \bbr$ given by $F_A(\x) =\langle A,
P_\x\rangle$.

\Prop{A.4}  {\sl Suppose $\f\in \L^p V^*$  and   $\s \in\End (V)$.
Then for all $p$-planes $\x\in G(p, V)$,}
$$
{\BM}_{\f} (\s)(\x)\ =\ (\tr_\x \s)\f(\x) + \f (\Der_{\wt A} \x)
\eqno(A.2)
$$
\pf
Pick an orthonormal basis for $\x$ and extend to an orthonormal basis of $V$.  
It suffices to prove (A.2) when   $A=b\otimes a$
with $a$ and $b$ elements of this basis. Using  formula (\ZZ.11) we see the following.

(1)\  If $a\in\x^\perp$, then all terms in (A.2) are zero.

(2)\  If $a\in\x$ and $b\in\x^\perp$, then $\wt A = A=b\otimes a$,   $\tr_\x A=0$,
and ${\BM}_{\f} (b\otimes a)(\x) = (a\wedge(b\hk\f))(\x) = \f(b\wedge(a\hk\x)) = \f(\Der_A\x)$

(3)\  If  $a=b\in\x$, then ${\BM}_{\f} (\s)(\x)= \f(a\wedge(a\hk\x)=\f(\x)$ and $\tr_\x(A)=1$.
Since $\wt A=0$, equation (A.2) holds in this case.

(4)\  If $a,b \in\x$ and $a\perp b$, then $b\wedge(a\hk \x)=0$, and one sees easily that
all three terms in (A.2) are zero.
\qed

 \Remark{A.5} Proposition A.2  can be restated as
 $$
 \lambda_{\f}^*(\x)\ =\ \f(\x) P_\x \qquad {\rm for\ all\ }\ \ \x\in G^{\rm cr}(\f).
 \eqno{(A.3)}
 $$
 Conversely, if $ \lambda_{\f}^*(\x)\ =\ c P_\x$\  for some  $\x\in G(p, V)$, then
 $c=\f(\x)$ and $\x$ is $\f$-critical.
 
 \pf  For all  $\x\in G(p,V)$ we have 
 $\langle  P_\x, \lambda_{\f}^*(\x)  \rangle =  ( \lambda_\f  P_\x)(\x)
  = (D_{P_{\x}^t}\f ) (\x)   =
  \f(D_{P_{\x}} \x)  = p\f(\x)$ 
since $D_{P_{\x}}\x=p\x$.  Therefore, $\lambda_{\f}^*(\x)=cP_\x$
 implies that $pc=p\f(\x)$ and equation (A.3) holds. 
 Equation (A.2) now implies that $\f(D_{\wt A}\x)=0$ for all $A\in \End(V)$ and, in particular,
  $\f(D_{L}\x)=0$ for all $L:\x\to \x^{\perp}$.  That is, $\f$ vanishes on $T_{\x}G\subset \L_p V$,
  i.e., $\x\in G^{\rm cr}(\f)$. \qed

\medskip

We now define an oriented submanifold $M$ of $X$ to be $\f$-{\bf critical} if $\oa T_x M \in G^{\rm cr}(\f)$ for all $x\in M$.  We leave it to the reader to use Proposition A.2  to establish the following extension of the previous results.

\Theorem{A.6} {\sl  Suppose $\f$ is a  $p$-form on a riemannian manifold 
$X$  and  $M\subset X$ is a $\f$-critical submanifold with mean curvature vector field $H$.     
Then for all $f\in C^\infty(X)$, 
$$
{\BM}_{\f} (\Hess f)\ =\ [\D_M (f) +H(f)]\f
$$
 when restricted to $M$. In particular, if $M$ is minimal, then on $M$}
 $$
 {\BM}_{\f} (\Hess f) \ =\ (\D_M f)\f
 $$
  
%
%
 %

\Ex{}  Let $\f={1\over 6}\{\o_I^2+\o_J^2+\o_K^2\}$ be the quaternion calibration on 
 $\bbh^n$. Then $\pm {1\over 3}$ are critical values and the  $\f$-critical submanifolds
 with critical value $\pm {1\over 3}$ include all complex Lagrangian submanifolds for any  complex
 structure defined by right multiplication by a unit imaginary quaternion (cf.  [U]).

 \Remark {A.7}  In a very interesting recent paper Colleen Robles [R] has shown that  for 
  any given parallel calibration $\phi$,  $\l_\f({\rm SkewEnd}(TX))$ generates an exterior differential 
  system whose integral submanifolds are exactly the $\phi$-critical submanifolds.
  At a point $x$ this can be stated equivalently as follows.
 Given $\x\in G_p(T_xX)$
  $$
  \l_\f(A)(\x) \ =\ 0 \ \ \forall\, A \in {\rm SkewEnd}(T_xX) \qquad  \iff \qquad \xi \in G^{\rm cr}(\f).
  \eqno{(A.4)}
  $$
This  can  be derived from (A.2) and (A.3).

   \vfill\eject

%
%

\centerline{\bf Appendix B.   Constructing $\phi$-plurisubharmonic  functions. }

\bigskip
  
 Straightforward calculation shows that if $F(x)=g(u_1(x),...,u_m(x))$,  then
 $$
\Hess F \ =\ \sum_{j=1}^m {\partial g \over \partial t_j} \Hess u_j
 + \sum_{i,j=1}^m {\partial^2 g \over \partial t_i\partial t_j} (\nabla u_i \circ \nabla u_j)
 \eqno {(B.1)}
 $$
 and hence
  $$
\ch^{\f}( F ) \ =\ \sum_{j=1}^m {\partial g \over \partial t_j} \ch^{\f}( u_j)
 + \sum_{i,j=1}^m {\partial^2 g \over \partial t_i\partial t_j} \BM_{\f}(\nabla u_i \circ \nabla u_j)
 \eqno {(B.1)'}
 $$

\Prop{B.1}  {\sl  If  $u_1,...,u_m$ are $\phi$-pluriharmonic  and $g(t_1,...,t_m)$ is convex,
then $F=g(u_1,...,u_m)$ is \fp. More generally, if ${\partial g \over \partial t_j} \geq 0$ for 
$j=1,...,m$ and $g$ is convex, then $F=g(u_1,...,u_m)$ is \fp whenever each $u_j$ is 
$\phi$-plurisubharmonic.
}
\pf Under our assumptions the  first term in equation (B.1)$'$ is $\geq 0$ on any $\xi\in G(\phi)$.
To show that the second term is $\geq 0$ is suffices to consider the case where the matrix
$(({\partial^2 g \over \partial t_i\partial t_j}))$ is rank one,
 i.e., equal to $((x_ix_j))$ for some vector $x\in\bbr^n$.
Then the second term equals $\BM_\phi\{(\sum_i x_i\nabla u_i)\circ (\sum_j x_j\nabla u_j)\}$
which is $\geq0$  on $\xi\in G(\phi)$ by (\ZZ.13) and Corollary \ZZ.7. \qed
\medskip

We now analyze the case where $m=2$ and determine necessary and sufficient 
conditions for $F=g(u_1,u_2)$ to be $\phi$-plurisubharmonic.

\Lemma{B.2} {\sl  Fix $v,w\in \bbr^n$ and $\xi\in G(\phi)$.  Let $v_0$ and $ w_0$ denote the orthogonal
projections of $v$ and $w$ respectively onto   $\xi$ (considered as a $p$-plane in $\bbr^n$). Then}
$$
\BM_\phi(v\circ w)(\xi)\ =\   \langle v_0, w_0\rangle  .
$$
\pf
Write $v=v_0+v_1$ and $w=w_0+w_1$ with respect to the decomposition $\bbr^n = \span \xi
\oplus (\span \xi)^\perp$.  Then for $\xi \in G(\phi)$ we have
$$
\BM_\phi(v\circ w)(\x) \ =\  \phi\{(v_0+v_1)\wedge ((w_0+w_1) \hk \xi))\}
 \ =\   \phi\{(v_0+v_1)\wedge (w_0 \hk \xi)\}
\qquad \qquad $$ $$\qquad 
  \ =\   \phi(v_0\wedge (w_0 \hk \xi))  \ =\    \langle v_0, w_0\rangle  
\phi(\xi)  \ =\  
 \langle v_0, w_0\rangle  .
$$
where the third equality follows from the First Cousin Principle. \qed \medskip

By Lemma B.2 we have that for  $\xi\in G(\phi)$, 
$$\eqalign{
\BM_\phi\{a\,v\circ v +&2b \, v\circ w +c \,w\circ w\} (\xi)
\cr
&=\ a\|v_0\|^2 +2b \langle v_0, w_0\rangle   + c\|w_0\|^2
\cr
&=  \  \left\langle \left( \matrix{a &b \cr  b & c}\right),   \left( \matrix{\|v_0\|^2 & \langle v_0, w_0\rangle   
 \cr    \langle v_0, w_0\rangle   & \|w_0\|^2 }\right) \right\rangle
}
\eqno{(B.2)}$$

\Remark{B.3} A symmetric $n\times n$-matrix $A$ is $\geq 0$ iff $ \langle A, P\rangle  \geq 0$ for all rank-one symmetric
$n\times n$-matrices $P$.

\Remark{B.4}  The matrix  $ \left( \matrix{\|v_0\|^2 & \langle v_0, w_0\rangle     \cr    \langle v_0, w_0\rangle  & \|w_0\|^2 }\right)$
is rank-one iff  $v_0$ and $w_0$ are linearly dependent.

\Lemma{B.5}  {\sl  Let $v,w\in \bbr^n$  be linearly independent. Suppose that for every line
$$
\ell\ \subset \ \span\{v,w\}
$$
there exists a $(p-1)$-plane $\xi_0\subset \span\{v,w\}^\perp$ such that
$\ell \oplus\xi_0$ (when properly oriented) is a $\phi$-plane.  Then $\BM_\phi\{a v\circ v+2b v\circ w +cw\circ w \}$ is $\phi$-positive if and only if  
$\left(\matrix{a &b\cr b&c}\right) \geq 0$.
}

\pf
Necessity is already done.  For sufficiency fix $a,b,c$.  For each $\ell \subset  \span\{v,w\}$
let $\xi\in G(\phi)$ be the oriented $p$-plane $\ell\oplus \xi_0$ given in the hypothesis, and
 note that by equation (B.2)
 $$
\BM_\phi\{a v\circ v+2b v\circ w +cw\circ w \} (\xi)\ =\  \left\langle \left( \matrix{a &b \cr  b & c}\right),                \left( \matrix{v_{\ell}^2 &v_{\ell}w_{\ell}  \cr   v_{\ell}w_{\ell} & w_{\ell}^2 }\right) \right\rangle
\ \geq\ 0
 $$
where $v_{\ell}= \langle v,e\rangle e$,  $w_{\ell}=  \langle w,e\rangle  e$, and $\ell = \span\{e\}$.  Now the matrix 
$\left( \matrix{v_{\ell}^2 &v_{\ell}w_{\ell}  \cr   v_{\ell}w_{\ell} & w_{\ell}^2 }\right) $
is rank-one, and every rank-one $2\times 2$ matrix, up to positive scalars, occurs in this family.
The result follows from Remark B.3. \qed

\Def{B.6}  A calibration $\phi$ on a manifold $X$ is called {\bf rich} (or {\bf 2-rich}) if for any 2-plane $P\subset  T_xX$ at any point x, and for any line $\ell\subset P$, there exists a $(p-1)$-plane
$\xi_0\subset P^\perp$ so that $\pm \ell\oplus \xi_0$ is a $\phi$-plane.

\Prop{B.7} {\sl   Let $(X,\phi)$ be a rich calibrated manifold.  Suppose   $u_1, u_2$ are 
$\phi$-pluriharmonic functions on $X$ 
with $\nabla u_1 \wedge \nabla u_2 \neq 0$ on a dense set. Then
for any $C^2$-function $g(t_1,t_2)$}
$$
F\ =\ g(u_1,u_2)  \in \fpsh \ \ \Leftrightarrow\ \ \ g {\rm \ \ is\ convex}
$$
\pf
Apply Proposition B.1, equation (B.2) and Lemma B.5.\qed

\Prop{B.8}  {\sl The Special Lagrangian calibration on a Calabi-Yau n-fold, $n\geq 3$, and the associative and coassociative calibrations on a $G_2$-manifold are rich calibrations.}

\pf
For the Special Lagrangian case it suffices to consider  $\phi ={\rm Re}(dz)$ on $\bbc^n$, $n\geq 3$.
Let $e_1, Je_1,..., e_n, Je_n$ be the standard hermitian basis of $\bbc^n$.  By unitary invariance we may assume that $\ell = \span\{e_1\}$ and $P=\span\{e_1, \a Je_1+ \b e_2\}$
(see [HL$_1$, Lemma 6.13] for example).
 Then the $(p-1)$-plane 
$\xi_0 = -Je_2\wedge Je_3\wedge e_4\wedge\cdots \wedge e_n$ does the job.

Consider now the associative calibration $\phi(x,y,z) = \langle x\cdot y, z\rangle  $ on the imaginary octonians
${\rm Im}({\bbo}) = {\rm Im}({\bbh}) \oplus \bbh\cdot \epsilon$ where $\bbh$ denotes the quaternions
and $\e$ is any unit vector in ${\rm Im}({\bbo})^\perp$.
By the transitivity of the group $G_2$ on $S^6=G_2/{\rm SU}(3)$ and the transitivity of SU(3)
on the tangent space, we may assume  $\ell = \span\{i\}$ and $P=\span\{i,j\}$ in Im$(\bbh)$.
We now choose $\xi_0 =  \epsilon\wedge ( i\cdot \epsilon)$.
For the coassociative calibration we choose $\xi_0 =k \wedge (i\epsilon)\wedge(k\epsilon)$
and note that $i\wedge \xi_0= i\wedge k \wedge (i\epsilon)\wedge(k\epsilon) $ is coassociative because its orthogonal complement   is $j\wedge\epsilon\wedge (j \epsilon) $ which is associative.

Alternatively, as noted by the referee, these latter cases follow easily from the fact that
G$_2$ acts transitively on the Stiefel manifolds $V_{2,7}$ of ordered pairs of
orthonormal 2-vectors in $\bbr^7$ (see for example [HL$_1$, Prop. IV.1.10]).
\qed

\medskip

We now give some examples and applications of the material above. We start with
Special Lagrangian geometry where the $\phi$-pluriharmonic functions are given by 
Proposition \AA.13.
Hence, we may apply Proposition B.7 to conclude the following. Let $u_1(z)$ and $u_2(z)$ be two
traceless hermitian quadratic forms on $\bbc^n$. (For example, take $u_1(z) =|z_1|^2-|z_2|^2$
and $u_2(z) =(n-2)|z_1|^2-|z_3|^2-\cdots-|z_n|^2$.)  Then $g(u_1(z),u_2(z))$ is \fp  if and only
if $g$ is convex.

\medskip
Formula (B.1)$'$ can be usefully applied to more general
functions $u_j$.   For example, in the Special Lagrangian case on
 $\bbc^n$ with $\phi={\rm Re}(dz)$, one has that $d d^\phi({1\over2}|z_k|^2) =\phi$, for any complex
 coordinate $z_k$  in any unitary coordinate system on $\bbc^n$. Hence  a linear combination of these functions has the property that $dd^\f u = c\f$ for some constant $c$.

\Prop{B.9} {\sl   Let $(X,\phi)$ be a rich calibrated manifold.  Suppose   $u_1, ..., u_n \in C^\infty(X)$  
 satisfy the equations $d d^\phi u_i= c_i \phi$ for constants $c_1,...,c_n$.  Then for any $C^2$-function $g(t_1,..., t_n)$}
$$
F\ =\ g(u_1,...,u_n)  \in \fpsh \ \ \Leftrightarrow\ \ \ \left\{\sum_{i=1}^nc_i{\partial g \over \partial t_i}\right\}
{\bf Id} \ +\   \left\langle \Hess_g ,  ((\langle (\nabla u_i)^\xi,  (\nabla u_j)^\xi  ))  \right\rangle             \ \geq \ 0
$$
for all $\phi$-planes $\xi$ at all points of $X$.
\vfill\eject


\centerline{\bf Appendix \B. Structure of the Core.}\medskip

Let $(X,\f)$ be a calibrated manifold and consider the set
$$
\cn \ \equiv\ \{ \x \in G(\f) : (\ch^\f f)(\x)\ =\ 0 \ \ {\rm for\ all\ }
f\in\fpsh\}.
$$
\Prop{\B.1} {\sl  Let $\pi:G(\f)\to X$ denote the projection. Then}
$$
\pi(\cn)\ =\ \Core(X).
$$
\pf
Suppose $x\notin\Core(X)$.  Then by definition there exists $f\in\fpsh$ 
with $(\ch^\f  f)(\x)>0$ for all $\x\in  \pi^{-1}(x)$.  Hence, $x\notin
\pi(\cn)$.

Conversely, suppose $x\notin \pi(\cn)$. Then for each $\x\in \pi^{-1}(x)$
there exists $f_\x\in\fpsh$ with $(\ch^\f f_\xi)(\x)>0$. Let $W_\x=\{\eta\in 
\pi^{-1}(x):  (\ch^\f f_\xi)(\eta)>0\}$ and choose a finite cover $W_{\x_1},..., 
W_{\x_\ell}$ of $\pi^{-1}(x)$. Then $f\equiv f_{\x_1}+\cdots+f_{\x_\ell}$
is strictly \fp at x, and so $x\notin \Core(X)$.  \qed

\Prop{\B.2} {\sl If $\x\in \cn$, then for each vector $v\in\span \x$,}
$$
df(v) \ =\ 0 \qquad{\sl  for\ all }\ \  f\in\fpsh
\eqno(\B.1)$$
\pf
Suppose $f\in\fpsh$ and set $F=e^f$.  Then $F\in \fpsh$, and   by  
 equation (B.1)' and Corollary \ZZ.7 we see that
$
0\ =\ (\ch^\f F)(\x)\ =\ e^f\{df\w d^\f f + \ch^\f f\}(\x)\ 
=\ e^f\{df\w d^\f f \} (\x) \ = \   
e^f|\n f\hk\x|^2
$.  \qed
\medskip

\Def{\B.3}  The {\bf tangential core} of $X$ is the set
$$
T\Core(X)\ \equiv\ \{ v\in TX: v \neq0\  {\rm and\ satisfies\ condition \
(\B.1)} \}. $$
Thus $T\Core(X) \subset TX$ is a subset defined by the vanishing of 
the family of smooth functions $df:TX\to \bbr$ for $f\in\fpsh$.
Propositions \B.1 and \B.2 show that the restriction of the bundle map
$p:TX\to X$ gives a surjective mapping
$$
p:T\Core(X) \to \Core(X)
$$
and for each $x\in X$, the vector space $T_x\Core(C)\equiv p^{-1}(x)$
contains the non-empty space generated by all $v\in \span \x$ for $\x\in
\cn_x$.

Consider a point $v\in T\Core(X)$ and suppose we have functions
$f_1,...,f_\ell\in\fpsh$ such that $\n df_1,..., \n df_\ell$ are linearly
independent at $v$. Then $T\Core(C)$ is locally contained in the
codimension-$\ell$ submanifold $\{df_1=\cdots=df_\ell=0\}$. 

\vfill\eject

\vskip .3in



\centerline{\bf References}

\vskip .2in

\noindent
[AFS]   B. S. Acharya, J. M. Figueroa-O'Farrill, B. Spence,  {\sl  Branes at angles and calibrated geometry},  JHEP 9804 (1998) 012.  
\smallskip

\noindent
[Al]   S. Alesker,  {\sl  Non-commutative linear algebra and  plurisubharmonic functions  of quaternionic variables}, Bull.  Sci.  Math., {\bf 127} (2003), 1-35.  

\smallskip

\noindent
[AV]   S. Alesker and M. Verbitsky,  {\sl  Plurisubharmonic functions  on hypercomplex manifolds and HKT-geometry}, J. Geom. Anal. {\bf 16} (2006), 375-399.

\smallskip

\noindent
[AF]  A. Andreotti and T. Frenkel,  {\sl  The Lefschetz theorem on hyperplane sections},  Ann. of Math. (2) {\bf 69}  (1959), 713-717.

\smallskip

\noindent
 [B]   R. L.  Bryant,  {\sl
 Calibrated cycles of codimension 3 in compact simple Lie groups},
  to  appear.
 \smallskip

\noindent
 [BH]   R. L.  Bryant and F. R. Harvey,  {\sl
 Submanifolds in hyperk\"ahler geometry},
 J.A.M.S., {\bf 2} (1989), 1-31.
 \smallskip

\noindent
 [BS]   R. L.  Bryant and S. M. Salamon,  {\sl
 On the construction of some complete metrics with exceptional holomony},
  Duke Math. J. {\bf 58} (1989),   829-850.

 \smallskip

\noindent
 [C]   E.  Calabi,  {\sl
 M\'etriques k\"ah\'eriennes et fibr\'es holomorphes},
 Annales scientifiques de l'\'Ecole Normale Superieure {\bf 12} (1979),   269-294.

 \smallskip

\noindent
 [DH]   J. Dadok and F. R. Harvey,  {\sl
Calibrations and spinors},
Acta Math.  {\bf 170} (1993),  83-119.

 \smallskip

\noindent
{[D]}  J.-P. Demailly, { Complex Analytic and Differential Geometry},
  e-book at Institut Fourier, UMR  5582 du CNRS,
   Universit\'e de Grenoble I, Saint-Martin d'H\`eres, France:
   can be found at http://www-fourier.ujfgrenoble.fr/~demailly/books.html.
 \smallskip

\noindent
 [DT]   S. K. Donaldson and R. P. Thomas,  {\sl  Gauge theory in higher dimensions}, in ``The Geometric Universe (Oxford 1996)'',  31-47, Oxford Univ. Press, Oxford, 1998.

 \smallskip

\noindent
 [EH]   T. Eguchi and A. J. Hanson,  {\sl
Asymptotically flat solutions to Euclidean gravity},
Physics Letters  {\bf 74B} (1978),   249-251.

 \smallskip

\noindent
 [EM]  ,   J. Evslin and L. Martucci, {\sl
D-brane networks in flux vacua, generalized cycles and calibrations},
JHEP0707:040, 2007. 	arXiv:hep-th/0703129v2.

 \smallskip

\noindent
 [Fu]   L. Fu, {\sl  On the boundaries of Special Lagrangian submanifolds},
Duke Math. J.   {\bf 79}   no. 2 (1995),   405-422.

 \smallskip

\noindent
[J]   D. D. Joyce,     
{    Compact Manifolds with Special Holonomy},    
Oxford University Press, Oxford, 2000.

\smallskip

\noindent
 [GL]  K. Galicki and B. Lawson, {\sl  Quaternionic reduction and quaternionic orbifolds}, Math. Ann. {\bf 282} (1989), 1-21.

 \smallskip

\noindent
[G]   J. P. Gauntlett,   {\sl Branes,  calibrations and supergravity}, 	arXiv:hep-th/0305074v3.

 \smallskip

\noindent
[GLW]   J. P. Gauntlett, N. D. Lambert and P. C. West,  {\sl Branes and calibrated geometries}, 
Commun. Math. Phys.  {\bf 202} (1999),  571-592.

 \smallskip

\noindent
[GP]   G.W. Gibbons and G. Papadopoulos,  {\sl Calibrations and intersecting branes}, 
Commun. Math. Phys.  {\bf 202} (1999), 593-619.

 \smallskip

\noindent
[GW]   F. Gmeiner and F. Witt,  {\sl Calibrated cycles and T-duality}, 
Commun. Math. Phys.  283:543-578, 2008. 	arXiv:math/0605710v5.

 \smallskip

\noindent
[GR]   H. Grauert and R. Remmert,  Coherent Analytic Sheaves, Springer-Verlag, Berlin-Heidelberg, 1984.
\smallskip

\noindent
[GZ]  V. Guedj and A. Zeriahi,     
{\sl    Intrinsic capacities on compact K\"ahler manifolds},    
J. Geom. Anal. {\bf 15} (2005), no. 4, 607-639.

\smallskip

\noindent
[H]  F.R. Harvey,
Spinors and Calibrations,  Perspectives in Mathematics, vol. 9 Academic Press, Boston, 1990

 \smallskip

   \noindent 
 {[HL$_1$]} F. R. Harvey and H. B. Lawson, Jr, {\sl Calibrated geometries},  Acta Mathematica 
{\bf 148} (1982), 47-157.

 \smallskip

   \noindent 
 {[HL$_2$]} F. R. Harvey and H. B. Lawson, Jr, {\sl Duality of positive currents and plurisubharmonic functions in calibrated geometry},  Amer. J. Math. (to appear).  ArXiv:math.0710.3921.
 \smallskip

   \noindent 
 {[HL$_3$]} F. R. Harvey and H. B. Lawson, Jr, {\sl Dirichlet duality and the nonlinear Dirichlet
 problem},  Comm.  Pure  Appl. Math.  {\bf  62} (2009),  396-443.
 \smallskip

   \noindent 
 {[HL$_4$]} F. R. Harvey and H. B. Lawson, Jr, {\sl   Plurisubharmonicity in a general geometric context},  (to appear).
 \smallskip

   \noindent 
 {[HL$_5$]} F. R. Harvey and H. B. Lawson, Jr, {\sl Lagrangian plurisubharmonicity and convexity 
 in Symplectic  geometry},  (to appear).
 \smallskip

   \noindent 
 {[HL$_6$]} F. R. Harvey and H. B. Lawson, Jr, {\sl Dirichlet duality and the nonlinear Dirichlet
 problem on riemannian manifolds},  (to appear).
 \smallskip



 \noindent
[HP] F. R. Harvey, J. Polking, {\sl Extending analytic objects},
Comm.  Pure Appl. Math. {\bf 28} (1975), 701-727.

 \smallskip

\noindent
[HW$_1$] F. R. Harvey,  R. O. Wells, Jr.,  {\sl Holomorphic approximation and hyperfunction 
theory on a $C^1$ totally real submanifold of a complex manifold},
  Math.  Ann. {\bf 197} (1972),  287-318.

 \smallskip

\noindent
[HW$_2$] F. R. Harvey,  R. O. Wells, Jr.,  {\sl Zero sets of non-negatively strictly plurisubharmonic
functions},
  Math.  Ann. {\bf 201} (1973),  165-170.

 \smallskip

\noindent
[Her] R. Hernandez,  {\sl Calibrated geometries and non perturbative superpotentials in M-theory},
 Eur. Phys. J. C18 (2001), 619-624.

 \smallskip

\noindent
{[Ho]}  L. H\"ormander,  {An introduction to the theory of functions of several complex variables},
     Van Nostrand Press, Princeton, N.J., 1966.

 \smallskip

   \noindent
[L$_1$]    H. B. Lawson, Jr.,      Minimal Varieties in Real and Complex Geometry, Les Presses de 
L'Universite de Montreal,  1974.

 \smallskip

 \noindent 
[L$_2$]   {\sl Minimal varieties}, Differential Geometry 
(Proc.  Sympos.  Pure Math., Vol XXVII, 
Stanford University, Stanford, Calif., 1973), Part 1, pp. 143-175.  Amer. Math. Soc. 
Providence, R. I., 1975.

 \smallskip

   \noindent
[MC]    O. Mac Conamhna,    {\sl  Spacetime singularity resolution by M-theory
fivebranes: calibrated geometry, Anti-de Sitter solutions and special holonomy metrics},    
Commun. Math. Phys. 284: 345-389, 2008. arXiv:0708.2568.
 \smallskip

   \noindent
[MS]    D. Martelli and J. Sparks,    {\sl  G-structures, fluxes and calibrations in M-theory},    
Phys. Rev. D68 (2003) 085014. 	arXiv:hep-th/0306225v2.
 \smallskip

 \noindent 
 {[M]} J. Milnor, {Morse Theory}, Annals of Math. Studies no. {\bf 51}, Princeton University Press,
 Princeton, N.J.,  1963.
 \smallskip

   \noindent
[OT]    Y. Ohnita and H. Tasaki,    {\sl  Uniqueness of certain 3-dimensional homologically volume minimizing submanifolds in compact simple Lie groups},    
Tsukuba J. Math. 
{\bf 10} (1986), 11-16.

 \smallskip

   \noindent
[ON]    B. O'Neill,    Semi-Riemannian Geometry,    Pure and Applied Math. no. 103, 
Academic Press, London, 1983.

 \smallskip

 \noindent 
 {[R]}  C. Robles, {A system of PDE for calibrated geometries}, ArXiv:0808.2158. \smallskip

   \noindent
[S]  H. H. Schaefer,  Topological Vector Spaces,    Springer Verlag,
New York,  1999.

\smallskip

   \noindent
[T]     H. Tasaki,    {\sl  Certain minimal or  homologically volume minimizing submanifolds in compact symmetric spaces},    
Tsukuba J. Math. 
{\bf 9} (1985), 117-131.

 \smallskip

   \noindent
[Th]     Dao \v Cong  Thi,    {\sl  Real minimal currents in compact Lie groups},    
Trudy Sem.  Vektor Tenzor. Anal. no. 19 (1979), 112-129.

 \smallskip

   \noindent
[Ti]   G. Tian,    {\sl  Gauge theory and calibrated geometry, I},    
Ann. of Math (2) {\bf 151} no. 1 (2000), 193-268.

 \smallskip

\noindent
[U]  I.  Unal, Ph.D. Thesis, Stony Brook, 2006.

   \noindent
[V]    M. Verbitsky,    {\sl  Manifolds with parallel differential forms and K\"ahler identities
for $G_2$-manifolds},  arXiv : math.DG/0502540 (2005).
 \smallskip

\end